\documentclass[10pt]{article}%

\usepackage[utf8]{inputenc}

\usepackage{amsmath}
\usepackage{amsfonts}
\usepackage{mathrsfs}
\usepackage{amssymb, color}
\usepackage[linkcolor=black,anchorcolor=black,citecolor=black]{hyperref}
\usepackage[numbers,sort&compress]{natbib}
\usepackage{graphicx}
\numberwithin{equation}{section}
\usepackage[body={15.5cm,21cm}, top=3cm]{geometry}%
\providecommand{\U}[1]{\protect\rule{.1in}{.1in}}
\providecommand{\U}[1]{\protect \rule{.1in}{.1in}}
\newtheorem{theorem}{Theorem}[section]

\newtheorem{definition}[theorem]{Definition}
\newtheorem{example}[theorem]{Example}

\newtheorem{lemma}[theorem]{Lemma}

\newtheorem{proposition}[theorem]{Proposition}
\newtheorem{remark}[theorem]{Remark}

\newtheorem{assumption}[theorem]{Assumption}
\newenvironment{proof}[1][Proof]{\noindent \textbf{#1.} }{\  \rule{0.5em}{0.5em}}

\def \E{\mathbf{E}}
\def \d{\mathsf{d}}

\def \P{\mathbf{P}}

\usepackage [numbers,sort&compress] {natbib} 

\def\d{\mathrm{d}}

\begin{document}
	\title{Mean Field Backward Stochastic Differential Equations with Double Mean Reflections}
	\author{ 	Hanwu Li \thanks{Research Center for Mathematics and Interdisciplinary Sciences, Shandong University, Qingdao 266237, Shandong, China. lihanwu@sdu.edu.cn.}
	\thanks{Frontiers Science Center for Nonlinear Expectations (Ministry of Education), Shandong University, Qingdao 266237, Shandong, China.}
    \and Jin Shi \thanks{Research Center for Mathematics and Interdisciplinary Sciences, Shandong University, Qingdao 266237, Shandong, China. }}
	\date{}
	\maketitle
	
	\begin{abstract}
     In this paper, we analyze the mean field backward stochastic differential equations ($\text{MFBSDEs}$) with double mean reflections, whose generator and constraints both depend on the distribution of the solution. When the generator is Lipschitz continuous, based on the backward Skorokhod problem with nonlinear constraints, we investigate the solvability of the doubly mean reflected $\text{MFBSDEs}$ by constructing a contraction mapping. Furthermore, if the constraints are linear, the solution can also be constructed by a penalization method. For the case of quadratic growth, we obtain the existence and uniqueness results by using a fixed-point argument, the BMO martingale theory and the $\theta$-method.
	\end{abstract}


    \textbf{MSC-classification}: 60H10
	
\section{Introduction}

In 1997, EI Karoui, Kapoudjian, Pardoux, Peng, and Quenez \cite{EKPPQ} first introduced the concept of reflected backward stochastic differential equations ($\text{RBSDEs}$),  characterized by the following dynamic equation: 
\begin{equation}\label{intro1}
 Y_t=\xi+\int_t^T f(s,Y_s,Z_s)ds-\int_t^T Z_s dB_s+(K_T-K_t),\quad \forall t\in[0,T].
\end{equation}
 $\text{RBSDEs}$  can be viewed as a variant of $\text{BSDEs}$,wherein the first component of the solution is constrained to remain above a specified process $L$, referred to as the obstacle, such that $Y_t\geq L_t$ for all $t\in[0,T]$. In contrast to the classical $\text{BSDEs}$, RBSDEs incorporate an additional nondecreasing predictable process $K$ which can be seen as the force aiming to push the solution upward. The primary aim is to identify the minimal solution $Y$, which is uniquely defined by the Skorokhod condition $\int_{0}^{T} (Y_t-L_t) dK_t=0$. Following this,  Cvitani\'{c} and Karatzas \cite{CK} explored BSDEs with two reflecting barriers, where the solution is required to remain within two specified continuous processes $L$ and $U$, designated as the lower and upper obstacles, respectively. Given its significance in both theoretical frameworks and practical applications, RBSDEs has garnered considerable scholarly interest (see, \cite{CM,GIOOQ,HL,K1,K2,KLQT,PX} and the associated references). 

In contrast to the previously discussed pointwise constraints on the solutions, Briand, Elie, and Hu  \cite{BEH} have recently investigated   $\text{BSDEs}$ with mean reflection, wherein the constraint is articulated in terms of the solution's distribution. Specifically, for a designated  loss function $l$, the constraint is formulated as $\mathbf{E}[l(t,Y_t)]\geq 0$, for any  $t\in[0,T]$. Unlike the classical reflected case, the force $K$ is required to be a deterministic function that adheres to the Skorokhod condition. The authors established the existence and uniqueness of what is termed the deterministic flat solution to mean reflected $\text{BSDEs}$, which they subsequently applied to address super-hedging problems within the context of ongoing risk management.  Following this, Falkowski and Słomiński \cite{FS} examined BSDEs with two reflecting constraints, focusing on mean reflection under dual constraints, specifically $$\mathbf{E}[l(t,Y_{t})]\in[l_t,r_t], \forall t\in[0,T].$$ In this scenario, $K$ is a deterministic function of bounded variation that satisfies the conditions
$$\int_{0}^{T}(\mathbf{E}[l(s,Y_s)]-l_s)dK_{s}^{+}=0, \int_{0}^{T}(\mathbf{E}[l(s,Y_s)]-r_s)dK_{s}^{-}=0$$
where $K^{+},K^{-}$ denote the positive and negative components of the Jordan decomposition of the function $K$. More recently, Li \cite{L} investigated a more general reflection constraint, that is, given two nonlinear loss functions $L,R$  with $L\leq R$, the first component $Y$ must satisfy the following condition $$\mathbf{E}[L(t,Y_t)]\leq 0\leq \mathbf{E}[R(t,Y_t)]$$ 
and the minimality conditions for the bounded variation function $K$ are defined accordingly.  For further insights into this area of research, we refer to works\cite{BH,GLX,HHLLW,LW,QF}.

A notable extension of BSDEs that incorporates the distribution of the solution is represented by the so-called $\text{MFBSDEs}$, as introduced by  Buckdahn, Djehiche, Li and Peng \cite{BDL}. This framework was developed in the context of analyzing the limiting behavior of a high-dimensional system of forward-backward BSDEs. Subsequently, Buckdahn, Li and Peng \cite{BLP} broadened the scope of their investigation into  $\text{MFBSDE}$ within a more general framework. They established key results including the existence and uniqueness of solution, the comparison theorem and the stochastic representation for a nonlocal $\text{PDE}$. For the case of reflected $\text{MFBSDEs}$, we may refer to \cite{Li',XZ} for the one-dimensional and  multi-dimensional cases, respectively. Both studies demonstrated the well-posedness of reflected $\text{MFBSDEs}$ and offered the probabilistic interpretation of the nonlocal PDEs with the obstacles. 

The present paper is dedicated to the study of the following $\text{MFBSDE}$ with double mean reflections over the time interval $[0,T]$ 
\begin{equation}\label{MFBSDEDMR}
\begin{cases}
Y_t=\xi+\int_t^T f(s,Y_s,\mathbf{P}_{Y_s},Z_s,\mathbf{P}_{Z_s})ds-\int_t^T Z_s dB_s+K_T-K_t, \\
\mathbf{E}[L(t,Y_t)]\leq 0\leq \mathbf{E}[R(t,Y_t)], \\
K_t=K^R_t-K^L_t,\ K^R,K^L\in I[0,T] \textrm{ and } \int_0^T \mathbf{E}[R(t,Y_t)]dK_t^R=\int_0^T \mathbf{E}[L(t,Y_t)]dK^L_t=0.
\end{cases}
\end{equation}
We will consider the well-posedness of Eq. \eqref{MFBSDEDMR} in various typical situations. In the context of Lipschitz generator and $p$-integrable terminal value, we demonstrate the uniqueness and existence by employing the relevant estimate of the Skorokhod problem and the contraction mapping principle. It is noteworthy that an alternative  effective approach for constructing the solution to a reflected $\text{BSDE}$ is the penalization method (refer to \cite{EKPPQ,CK,PX}). Specifically, in the case of mean reflection with a linear constraint, the solution can also be approximated by a sequence of penalized $\text{MFBSDEs}$. 

Given the complexities associated with $\text{MFBSDEs}$ with quadratic generator depending on the distribution of $Z$, this study will focus on $\text{MFBSDEs}$ with mean reflections, specifically considering the quadratic generator that relies solely on the distribution of the process $Y$. The foundational element for constructing the solution is expressed as  $Y_{t}=y_{t}+K_{T}-K_{t}$, where $y$ is the first component of an appropriate BSDE and $K$ denotes the solution to a specific backward Skorokhod problem. By employing estimate related to the backward Skorokhod problem, the \text{BMO} martingale theory and the contraction mapping principle, we demonstrate the well-posedness of doubly mean reflected $\text{MFBSDEs}$ for the case of quadratic generator and bounded terminal value. When the terminal value is unbounded, we draw upon insights from references \cite{BH08} and \cite{HMW} to introduce a technical Lemma \ref{le-2}, which plays a crucial role in the proof of both uniqueness and existence for Eq. \eqref{MFBSDEDMR}. Based on this lemma, we apply estimate pertinent to Skorokhod problem and the $\theta$-method to establish the uniqueness and the solution is constructed by the Picard iteration method. However, the sequence of $\left\{K^{(m)}\right\}_{m=0}^{\infty}$ we constructed by Picard iteration is a family of bounded variation functions and proving the existence of the solution for Eq. \eqref{MFBSDEDMR} requires a uniformly estimate for $K^{(m)}$. We introduce the new backward Skorokhod problem $\mathbb{BSP}^{r^{1,(m)}}_{l^{1,(m)}}(s^{1,(m)},a^{1,(m)})$ with two linear constraints. Combine with Proposition $\ref{totalvariation}$ and $\ref{comparison}$ in this paper to obtain the desired conclusion. In addition, the challenges arise in verifying that the limit function $K$ adheres to the Skorokhod conditions, as it is not monotonic, as noted in \cite{HMW}. To address this issue, it is sufficient to demonstrate that $K$ is equivalent to $K'$, where $K'$ is the final component of the solution to doubly mean reflected BSDE with a constant generator $f(s,Y_s,\mathbf{P}_{Y_s},Z_s)$ and loss functions $L,R$. The equivalence can be established through estimate derived from the backward Skorokhod problem.

Actually, the reflected BSDEs, where both the generator and the constraint are contingent upon the distribution of the solution, have been investigated in the literature by \cite{DEH,HMW',HMW} for the single reflection case and by \cite{CHM} for the double reflections case. In the work presented in \cite{CHM}, the constraints are defined as follows
\begin{align*}
h(t,\omega, Y_t, \mathbf{E}[Y_t])\leq Y_t \leq g(t,\omega, Y_t, \mathbf{E}[Y_t]),
\end{align*}
which cannot infer the constraints made in this paper. It is noteworthy  that the generator discussed in \cite{CHM} does not depend on the distribution of the $Z$-term and is confined to the Lipschitz continuous framework. Additionally, in order to demonstrate the existence and uniqueness of solution, an extra assumption concerning the Lipschitz constants is required.  The penalty method is also explored in \cite{CHM}, which needs nondecreasing property of the generator. In contrast, the present paper does not impose these supplementary conditions. 

This paper is organized as follows. In Section 2, we recall some basic results about the Skorokhod problem and backward Skorokhod problem. In Section 3,  we study the existence and uniqueness of doubly mean reflected MFBSDEs with the Lipschitz generator. In Section 4, doubly mean reflected MFBSDEs with quadratic generator are investigated for both bounded and unbounded terminal values. The construction via penalization for the Lipschitz generator and some technical proofs are put into the appendix.  
 
\section{Preliminaries}

 Fix a finite time horizon $T>0$. Consider a filtered probability space $(\Omega,\mathcal{F},\mathbb{F},\P)$ satisfying the usual conditions of right continuity and completeness. Let $B$ be a $1$-dimensional standard Brownian motion. Actually, all the results in this paper still hold when we consider the multi-dimensional Brownian motion. For each $p>1$ and for any finite time interval $[u,v]\subset [0,T]$, we first introduce the following notations.
 
\begin{itemize}	
\item $\mathcal{L}^p_v$: the collection of real-valued $\mathcal{F}_{v}$-measurable random variable $\xi$ satisfying 
	$$\left \| \xi \right \|_{\mathcal{L}^p}=\mathbf{E} \left [ \left | \xi  \right | ^p \right ]^{\frac{1}{p} }< \infty. $$
	 
	\item $\mathcal{L}^{\infty}_v$: the collection of real-valued $\mathcal{F}_{v}$-measurable random variable $\xi$ satisfying
	 $$\left \| \xi  \right \|_{\mathcal{L}^\infty}=\operatorname{ess\,sup}\left |\mathbf{\xi} \left ( \omega  \right )   \right |<\infty. $$
	
	\item $\mathcal{H}^{p}\left[u,v\right]$: the collection of real-valued $\mathcal{F}$-progressively measurable processes $\left(z_{t}\right)_{u\leq t\leq v}$ satisfying 
	$$\left \| z \right \|_{\mathcal{H}^p }=\mathbf{E}\left [\left ( \int_{u}^{v}\left | z_{t} \right |^2 d t \right )^{\frac{p}{2}}   \right ]^{\frac{1}{p} }<\infty. $$
	
	\item $S^p\left[u,v\right]$: the collection of real-valued $\mathcal{F}$-adapted continuous processes $\left(y_{t}\right)_{u\leq t \leq v}$ satisfying 
	$$ \left \| y \right \|_{\mathcal{S}^p }=\mathbf{E}\left [ \sup _{t\in[u,v]} \left | y_{t} \right |^p \right ]^{\frac{1}{p}}  <\infty.$$
	
	\item $S^{\infty}\left[u,v\right]$: the collection of real-valued $\mathcal{F}$-adapted continuous processes $\left(y_{t}\right)_{u\leq t \leq v}$ satisfying
	$$\left \| y \right \|_{\mathcal{S}^{\infty } } =\operatorname{ess\,sup}_{(t,\omega)\in\left[u,v\right]\times\Omega}\left | y_t\left (\omega\right )  \right |  <\infty.$$
	
	\item $\mathcal{P}_{p}\left(\mathbb{R}\right)$: the collection of all probability measures over $\left(\mathbb{R},\mathcal{B}\left(\mathbb{R}\right)\right)$ with finite $p$-th moment, endowed with the $p$-Wasserstein distance $W_{p}$, where $p > 1$ is a positive integer and $$W_{p}(\mu,\nu)=\left(\inf_{X \sim \mu,Y \sim\nu}\mathbf{E}[\left|X-Y\right|]^p\right)^{\frac{1}{p}}.$$
	
	\item $\mathbb{L}^p_v$: the collection of real-valued $\mathcal{F}_{v}$-measurable random variable $\xi$ satisfying $\mathbf{E}\left[e^{p\left|\xi\right|}\right]<\infty.$
	
	\item $\mathbb{S}^p[u,v]$: the collection of all stochastic processes $Y$ such that $e^Y\in\mathcal{S}^p\left[u,v\right].$
	
	\item $\mathbb{L}_v$: the collection of all random variable $\xi\in\mathbb{L}_v^p$ for any $p> 1$, and $\mathcal{H}\left[u,v\right]$ and $\mathbb{S}\left[u,v\right]$ are defined similarly.
	
	\item $\mathcal{T}_{t}[u,v]$: the collection of $\left[u,v\right]$-valued $\mathcal{F}$-stopping times $\tau$ such that $\tau \geq t$, $\mathbf{P}$-a.s. 
	
	\item $BMO\left[u,v\right]$: the collection of real-valued progressively measurable processes $\left(z_{t}\right)_{u\leq t \leq v}$ such that 
	$$\left \| z \right \|_{BMO}:=\sup_{\tau\in \mathcal{T}_{u}}\operatorname{ess\,sup}_{\omega\in\Omega}\mathbf{E}_{\tau}\left[\int_{\tau}^{v}\left | z_{s} \right |^2 d s \right]^{\frac{1}{2}}<\infty. $$
	
	\item $C\left[u,v\right]$: the set of continuous functions from $\left[u,v\right]$ to $\mathbb{R}$.

 
	\item  $BV\left[u,v\right]$: the set of functions $K\in C\left[u,v\right]$ with $K_u=0$ and $K$ is of bounded variation on $\left[u,v\right]$.
	
	\item  $I\left[u,v\right]$: the set of functions in $C\left[u,v\right]$ starting from the origin which is nondecreasing.
 \end{itemize}
	
When the interval $[u,v]=[0,T]$, we always omit the time index. 
 For example, we write $\mathcal{H}^{p},\mathcal{S}^p,\mathcal{S}^{\infty}$ and so on. 

\subsection{Skorokhod problem}
\begin{definition}\label{def1}
Let $s\in C[0,T]$,  and $l,r:[0,T] \times \mathbb{R}\rightarrow \mathbb{R}$ be two functions with $l\leq r$.  A pair of functions $(x,K)\in C[0,T]\times BV[0,T]$ is called a solution to the  Skorokhod problem for $s$ with nonlinear constraints $l,r$ ($(x,K)=\mathbb{SP}_l^r(s)$ for short) if 
\begin{itemize}
\item[(i)] $x_t=s_t+K_t$;
\item[(ii)] $l(t,x_t)\leq 0\leq r(t,x_t)$;
\item[(iii)]  $K_{0-}=0$ and $K$ has the decomposition $K=K^r-K^l$, where $K^r,K^l$ are nondecreasing functions satisfying  
\begin{align}
\int_0^{T} I_{\{l(s,x_s)<0\}}dK^l_s=0, \  \int_0^{T} I_{\{r(s,x_s)>0\}}dK^r_s=0.
\end{align}
\end{itemize}
\end{definition}

We propose the following assumption on the functions $l,r$.
\begin{assumption}\label{asslr}
The functions $l,r:[0,T] \times \mathbb{R}\rightarrow \mathbb{R}$ satisfy the following conditions
\begin{itemize}
\item[(i)] For each fixed $x\in\mathbb{R}$, $l(\cdot,x),r(\cdot,x)\in C[0,T]$;
\item[(ii)] For any fixed $t\geq 0$, $l(t,\cdot)$, $r(t,\cdot)$ are strictly increasing;
\item[(iii)] There exist two positive constants $0<c<C<\infty$, such that for any $t\in [0,T]$ and $x,y\in \mathbb{R}$,
\begin{align*}
&c|x-y|\leq |l(t,x)-l(t,y)|\leq C|x-y|,\\
&c|x-y|\leq |r(t,x)-r(t,y)|\leq C|x-y|.
\end{align*} 
\item[(iv)] $\inf_{(t,x)\in[0,T]\times\mathbb{R}}(r(t,x)-l(t,x))>0$.
\end{itemize}
\end{assumption}

\begin{theorem}[\cite{Li}]\label{SP}
Suppose that $l,r$ satisfy Assumption \ref{asslr}. For any given $s\in C[0,T]$, there exists a unique pair of solution to the Skorokhod problem $(x,K)=\mathbb{SP}_l^r(s)$.
\end{theorem}

The following proposition provides the continuity property of the solution  to the Skorokhod problem with respect to input function $s$ and the reflecting boundary functions $l,r$.

\begin{proposition}[\cite{Li}]\label{continuity}
Suppose that $(l^i,r^i)$ satisfy Assumption $\ref{asslr}$, $i=1,2$.
Given  $s^i\in C[0,T]$, let $(x^i,K^i)$ be the solution to the Skorokhod problem $\mathbb{SP}_{l^i}^{r^i}(s^i)$. Then, we have
\begin{equation}\label{diffk}
\sup_{t\in[0,T]}\left|K^1_t-K^2_t\right|
\leq \frac{C}{c}\sup_{t\in[0,T]}\left|s^1_t-s^2_t\right|+\frac{1}{c}(\bar{L}_T\vee\bar{R}_T),
\end{equation}
where
\begin{align*}
\bar{L}_T=\sup_{(t,x)\in[0,T]\times \mathbb{R}}\left|{l}^1(t,x)-{l}^2(t,x)\right|,\
\bar{R}_T=\sup_{(t,x)\in[0,T]\times \mathbb{R}}\left|{r}^1(t,x)-{r}^2(t,x)\right|.
\end{align*}
\end{proposition}


In the following, we aim to provide an estimate for the total variation of $K$, where $K$ is the second component of the solution to the Skorokhod problem $\mathbb{SP}_l^r(s)$. To this end, let $\phi_t,\psi_t$ be the unique solution to the following equations, respectively 
$$l(t,s_t+x)=0,\ r(t,s_t+x)=0.$$
Use $Var^{b}_{a}(f)$ to represent the total variation of function $f$ on time interval $[a,b]$.
By Theorem 2.2 and Theorem 2.3 in \cite{Li}, we have the following result. 
\begin{proposition}\label{totalvariation}
    Under the same assumption as in Theorem \ref{SP}, let $(x,K)$ be the solution to the Skorokhod problem $\mathbb{SP}_l^r(s)$. Then, for any $T>0$, we have
    $$Var^{T}_{0}(K)=K^{r}_{T}+K^{l}_{T}\leq Var^{T}_{0}(\phi)+Var^{T}_{0}(\psi).$$  
\end{proposition}

Next, we give the comparison property of solutions to Skorokhod problems with respect to the reflecting boundary functions.
\begin{proposition}[\cite{Li}]\label{comparison}
Suppose  $(l^i,r^i)$ satisfy Assumption $\ref{asslr}$ for $i=1,2$ with $r^1\geq r^2$ and $l^1\leq l^2$. For any given $s\in C[0,T]$, let $(x^i,K^i)$ solve the Skorokhod problem $\mathbb{SP}_{r^i}^{l^i}(s)$ with $K^i=K^{i,r}-K^{i,l}$. Then, for any $t\geq 0$, 
$$K^{2,r}_t\geq K^{1,r}_t \text{and} \ K^{2,l}_t\geq K^{1,l}_t. $$
\end{proposition}

\begin{remark}
    Actually, Theorem \ref{SP} and Proposition \ref{continuity} still hold if the input function and the reflecting boundary functions are c\`{a}dl\`{a}g. The main difference is that each component of the solution $(x,K)$ to the Skorokhod problem is c\`{a}dl\`{a}g.
\end{remark}

\subsection{Backward Skorokhod problem}

\begin{definition}\label{def2}
Let $s\in C[0,T]$, $a\in \mathbb{R}$ and $l,r:[0,T]\times \mathbb{R}\rightarrow \mathbb{R}$ be two functions such that $l\leq r$ and $l(T,a)\leq 0\leq r(T,a)$. A pair of functions $(x,K)\in C[0,T]\times BV[0,T]$ is called a solution of the backward Skorokhod problem for $s$ with nonlinear constraints $l,r$ ($(x,K)=\mathbb{BSP}_l^r(s,a)$ for short) if 
\begin{itemize}
\item[(i)] $x_t=a+s_T-s_t+K_T-K_t$;
\item[(ii)] $l(t,x_t)\leq 0\leq r(t,x_t)$, $t\in[0,T]$;
\item[(iii)]  $K$ has the decomposition $K=K^r-K^l$, where $K^r,K^l\in I[0,T]$ satisfy 
\begin{align}\label{iii}
\int_0^T I_{\{l(s,x_s)<0\}}dK^l_s=0, \  \int_0^T I_{\{r(s,x_s)>0\}}dK^r_s=0.
\end{align}
\end{itemize}
\end{definition}


\begin{theorem}[\cite{L}]\label{BSP}
Let Assumption \ref{asslr} holds. For any given $s\in C[0,T]$ and $a\in \mathbb{R}$ with $l(T,a)\leq 0\leq r(T,a)$, there exists a unique solution to the backward Skorokhod problem $(x,k)=\mathbb{BSP}_l^r(s,a)$. 
\end{theorem}
\begin{remark}\label{BSPtoSP}
Set
\begin{equation}
    \begin{split}
     &\bar{s}_t=a+s_T-s_{T-t},\ x_t=\bar{x}_{T-t},\ K_t=\bar{K}_{T}-\bar{K}_{T-t}, \ t\in[0,T],\\ 
     &\bar{l}\left(t,x\right)=l\left(T-t,x\right),\ \bar{r}\left(t,x\right)=r\left(T-t,x\right), \ (t,x)\in[0,T]\times\mathbb{R}.
    \end{split}
\end{equation}
According to Theorem 3.9 in \cite{L},  $\left(\bar{x},\bar{K}\right)$ is the unique solution to the Skorokhod problem $\mathbb{SP}^{\bar{r}}_{\bar{l}}\left(\bar{s}\right)$.
\end{remark}   
 
The following proposition provides the continuous dependence of the solution with respect to the input function $s$ and the reflecting boundary functions $l,r$.

\begin{proposition}[\cite{L}]\label{continuous}
Given $a^i\in\mathbb{R}$, $s^i\in C[0,T]$, $l^i,r^i$ satisfy Assumption \ref{asslr} and $l^i(T,a^i)\leq 0\leq r^i(T,a^i)$, $i=1,2$, 
 let $(x^i,k^i)$ be the solution to the backward Skorokhod problem $\mathbb{BSP}_{l^i}^{r^i}(s^i,a^i)$, $i=1,2$. Then, we have
\begin{equation}\label{diffK}
\sup_{t\in[0,T]}\left|K^1_t-K^2_t\right|
\leq  2\frac{C}{c}\left|a^1-a^2\right|+4\frac{C}{c}\sup_{t\in[0,T]}\left|s^1_t-s^2_t\right|+\frac{2}{c}(\bar{L}_T\vee\bar{R}_T),
\end{equation}
where $\bar{L}_T$ and $\bar{R}_T$ are the same as in Proposition \ref{continuity}.
\end{proposition}

\section{Doubly mean reflected MFBSDEs with Lipschitz generator}

The main purpose of this paper is to study the MFBSDE with double mean reflections of the following type
\begin{equation}\label{nonlinearyz}
\begin{cases}
Y_t=\xi+\int_t^T f(s,Y_s,\P_{Y_s},Z_s,\P_{Z_s})ds-\int_t^T Z_s dB_s+K_T-K_t, \\
\E[L(t,Y_t)]\leq 0\leq \E[R(t,Y_t)], \\
K_t=K^R_t-K^L_t,\ K^R,K^L\in I[0,T] \textrm{ and } \int_0^T \E[R(t,Y_t)]dK_t^R=\int_0^T \E[L(t,Y_t)]dK^L_t=0.
\end{cases}
\end{equation}

Compared with the doubly mean reflected BSDE studied in \cite{L}, the generator $f$  of the MFBSDE with double mean reflections depends on the distribution of the solution, which is required to satisfy the following assumption.
\begin{assumption}\label{ass1} 
    The generator $f$ is a map from $\Omega\times[0,T]\times \mathbb{R}\times \mathcal{P}_1(\mathbb{R})\times\mathbb{R}\times \mathcal{P}_1(\mathbb{R})$ to $\mathbb{R}$. For 
each fixed $(y,\mu,z,\nu)\in  \mathbb{R}\times \mathcal{P}_1(\mathbb{R})\times\mathbb{R}\times \mathcal{P}_1(\mathbb{R})$, $f(\cdot,\cdot,y,\mu,z,\nu)\in \mathcal{H}^p$ with some $p>1$. There exists $\lambda>0$ such that for any $t\in[0,T]$ and any $y,y',z,z'\in\mathbb{R}$, $\mu,\mu',\nu,\nu'\in\mathcal{P}_1(\mathbb{R})$
\begin{align*}
\left|f(t,y,\mu,z,\nu)-f(t,y',\mu',z',\nu')\right|\leq \lambda\left(\left|y-y'\right|+\left|z-z'\right|+W_1(\mu,\mu')+W_1(\nu,\nu')\right).
\end{align*}
\end{assumption}

The loss functions $L,R:\Omega\times [0,T]\times\mathbb{R}\rightarrow \mathbb{R}$ are measurable maps with respect to $\mathcal{F}_T\times \mathcal{B}([0,T])\times \mathcal{B}(\mathbb{R})$ satisfying the following conditions.

\begin{assumption}\label{ass2}
\begin{itemize}
\item[(1)] For any fixed $(\omega,x)\in \Omega\times\mathbb{R}$, $L(\omega, \cdot, x), R(\omega, \cdot, x)$ are continuous;
\item[(2)]  $\E\left[\sup_{t\in[0,T]}\left|L(t,0)\right|\right]<\infty$, $\E\left[\sup_{t\in[0,T]}\left|R(t,0)\right|\right]<\infty$;
\item[(3)] For any fixed $(\omega,t)\in \Omega\times [0,T]$, $L(\omega,t,\cdot),R(\omega,t,\cdot)$ are strictly increasing and there  exist two constants $c,C$ satisfying $0<c<C$ such that for any $x,y\in \mathbb{R}$,
\begin{align*}
&c\left|x-y\right|\leq \left|L(\omega,t,x)-L(\omega,t,y)\right|\leq C\left|x-y\right|,\\
&c\left|x-y\right|\leq \left|R(\omega,t,x)-R(\omega,t,y)\right|\leq C\left|x-y\right|;
\end{align*}
\item[(4)] $\inf_{\omega,t,x} \left(R(\omega,t,x)-L(\omega,t,x)\right)>0$.
\end{itemize}
\end{assumption}

\begin{remark}
 If $L\equiv -\infty$, the $\text{MFBSDE}$ with double mean reflections degenerates into the general mean reflected BSDE studied in \cite{HMW}.

\end{remark}

Before investigating the doubly mean reflected MFBSDE \eqref{nonlinearyz}, we first study the case of constant generator. 
\begin{proposition}\label{pro-1}
   Suppose that $\xi\in\mathcal{L}_T^p$ satisfies $\E[L(T,\xi)]\leq 0 \leq \E[R(T,\xi)]$ and Assumption \ref{ass2} holds. Given $C\in \mathcal{H}^p$, the $\text{BSDE}$ with double mean reflections 
   \begin{equation}\label{MFBSDEC}
   \begin{cases}
Y_t=\xi+\int_t^T C_s ds-\int_t^T Z_s dB_s+K_T-K_t, \\
\E[L(t,Y_t)]\leq 0\leq \E[R(t,Y_t)], \\
K_t=K^{R}_t-K^{L}_t, \int_0^T \E[R(t,Y_t)]dK_t^{R}=\int_0^T \E[L(t,Y_t)]dK^{L}_t=0,
\end{cases}
\end{equation}
has a unique solution $\left(Y,Z,K\right)\in \mathcal{S}^p\times\mathcal{H}^{p}\times BV[0,T]$.
\end{proposition}
\begin{proof}
 The proof is similar to the one for Proposition 4.2 in \cite{L}. For the purpose of the remaining proof, we briefly describe the construction here.
 
Let $(\tilde{Y},Z)\in \mathcal{S}^p\times \mathcal{H}^p$ be the solution to the BSDE with terminal value $\xi$ and constant generator $C$. For any $t\in[0,T]$, set 
 $$\begin{aligned}
  & s_t=\E\left[\int_{0}^{t}C_s ds\right], a=\E\left[\xi\right]
 \end{aligned}$$
and for any $(t,x)\in[0,T]\times \mathbb{R}$, we define 
$$l(t,x):=\E\left[L\left(t,\tilde{Y}_t-\E[\tilde{Y}_t]+x\right)\right],r(t,x):=\E\left[R\left(t,\tilde{Y}_t-\E[\tilde{Y}_t]+x\right)\right].$$
It is easy to check that $s$ is a continuous function  and $l,r$ satisfy Assumption \ref{asslr}. By Theorem \ref{BSP}, the backward Skorokhod problem $\mathbb{BSP}_l^r(s,a)$ admits a unique solution $(x,K)$. Set$$Y_t=\tilde{Y}_t+K_T-K_t=\xi+\int_{t}^{T}C_s ds-\int_{t}^{T}Z_s d B_s+K_T-K_t.$$ 
Then, $(Y,Z,K)$ is the solution to \eqref{MFBSDEC}.
\end{proof}

We are now ready to state the main result of this section.
\begin{theorem}\label{th-1}
Suppose  that $\xi\in\mathcal{L}_T^p$ satisfies $\E[L(T,\xi)]\leq 0 \leq \E[R(T,\xi)]$ and Assumptions $\ref{ass1}$, $\ref{ass2}$ hold. Then,  the $\text{MFBSDE}$  with double mean reflections (\ref{nonlinearyz}) admits a unique solution $(Y,Z,K)\in \mathcal{S}^p \times \mathcal{H}^{p} \times BV[0,T]$.
\end{theorem}
\begin{proof}
In this proof,  $M$ always represents a positive constant depending on $c,C,\lambda,p$, which may vary from line to line. Given $U^i\in \mathcal{S}^p$, $V^i\in \mathcal{H}^{p}$, $i=1,2$, Proposition $\ref{pro-1}$ ensures that the following $\text{BSDE}$ with double mean reflections admits a unique solution $(Y^i,Z^i,K^i)\in \mathcal{S}^p \times \mathcal{H}^{p} \times BV[0,T]$,
\begin{displaymath}
\begin{cases}
Y^i_t=\xi+\int_t^T f(s,U^i_s,\P_{U^i_s},V^i_s,\P_{V^i_s})ds-\int_t^T Z^i_s dB_s+K^i_T-K^i_t, \\
\E[L(t,Y^i_t)]\leq 0\leq \E[R(t,Y^i_t)], \\
K^i_t=K^{i,R}_t-K^{i,L}_t, \int_0^T \E[R(t,Y^i_t)]dK_t^{i,R}=\int_0^T \E[L(t,Y^i_t)]dK^{i,L}_t=0.
\end{cases}
\end{displaymath}
We define
\begin{align*}
\Gamma:\mathcal{S}^p\times\mathcal{H}^{p}&\rightarrow \mathcal{S}^p\times\mathcal{H}^{p},  \\
\Gamma(U^i,V^i)&=(Y^i,Z^i).
\end{align*}
Set
\begin{align*}
\hat{F}_t=F^1_t-F^2_t, \textrm{ where } F=Y,Z,K,U,V, \ \hat{f}_t=f(t,U^1_t,\P_{U^1_t},V^1_t,\P_{V^1_t})-f(t,U^2_t,\P_{U^2_t},V^2_t,\P_{V^2_t}). 
\end{align*}
Using the following representation for $1$-Wasserstein distance
\begin{align*}
W_1(\mu,\nu)=\inf_{X \sim\mu, Y\sim \nu}\E[\left|X-Y\right|],
\end{align*}
it is easy to check that 
\begin{align*}
|\hat{Y}_t|&=\left|\E_t\left[\int_t^T \hat{f}_sds\right]+\hat{K}_T-\hat{K}_t\right|\\
&\leq \lambda \E_t\left[\int_0^T \left(|\hat{U}_s|+|\hat{V}_s|+W_1\left(\P_{U^1_s},\P_{U^2_s}\right)+W_1\left(\P_{V^1_s},\P_{V^2_s}\right)\right)ds \right]+|\hat{K}_T|+|\hat{K}_t|\\
&\leq \lambda \E_t\left[\int_0^T \left(|\hat{U}_s|+|\hat{V}_s|+\E[|\hat{U}_s|]+\E[|\hat{V}_s|]\right)ds \right]+|\hat{K}_T|+|\hat{K}_t|.
\end{align*}
Applying the Doob's maximal inequality, we obtain that
\begin{equation}\begin{split}\label{differY}
 \E[\sup_{t\in[0,T]}|\hat{Y}_t|^p]&\leq M \E\left[\sup_{t\in[0,T]}\left(\E_{t}\left[\int_{0}^{T}|\hat{U}_s|+|\hat{V}_s|+\E[|\hat{U}_s|]+\E[|\hat{V}_s|]ds\right]\right)^p\right]+M \sup_{t\in[0,T]}|\hat{K}_t|^p\\
 &\leq M \E\left[\left(\int_{0}^{T}\left(|\hat{U}_s|+|\hat{V}_s|+\E[|\hat{U}_s|]+\E[|\hat{V}_s|]\right)ds\right)^p\right]+M \sup_{t\in[0,T]}|\hat{K}_t|^p.\\   
\end{split}
\end{equation}
Recalling the proof of Proposition \ref{pro-1} and \eqref{diffK}, we have
\begin{align*}
\sup_{t\in[0,T]}|\hat{K}_t|\leq  M \left\{\sup_{t\in[0,T]}|\hat{s}_t|+\sup_{(t,x)\in[0,T]\times\mathbb{R}}|\hat{l}(t,x)|\vee \sup_{(t,x)\in[0,T]\times\mathbb{R}} |\hat{r}(t,x)|\right\},
\end{align*}
where $\hat{s}_t=s^1_t-s^2_t$, $\hat{l}(t,x)=l^1(t,x)-l^2(t,x)$, $\hat{r}(t,x)=r^1(t,x)-r^2(t,x)$ and for $i=1,2$, 
\begin{align*}
&s^i_t=\E\left[\int_0^t f\left(s,U^i_s,\P_{U^i_s},V^i_s,\P_{V^i_s}\right)ds\right], \\ 
&l^i(t,x)=\E\left[L\left(t,\widetilde{Y}^i_t-\E[\widetilde{Y}^i_t]+x\right)\right],\\  &r^i(t,x)=\E\left[R\left(t,\widetilde{Y}^i_t-\E[\widetilde{Y}^i_t]+x\right)\right].
\end{align*}
Here, $\widetilde{Y}^i$ is the first component of the solution to the $\text{BSDE}$ with terminal value $\xi$ and constant driver $\left\{f\left(s,U_s^i,\P_{U_s^i},V_s^i,\P_{V^i_s}\right)\right\}_{s\in[0,T]}$, $i=1,2$. By Assumptions \ref{ass1} and \ref{ass2}, noting that $\widetilde{Y}^i_t=\E_t\left[\xi+\int_t^T f\left(s,U_s^i,\P_{U_s^i},V_s^i,\P_{V^i_s}\right)d s\right]$, we obtain that 
\begin{align*}
\left|l^1(t,x)-l^2(t,x)\right|&\leq C\E\left[\left|\left(\widetilde{Y}^1_t-\E[\widetilde{Y}^1_t]\right)-\left(\widetilde{Y}^2_t-\E[\widetilde{Y}^2_t]\right)\right|\right]\leq 2C\E[|\widetilde{Y}^1_t-\widetilde{Y}^2_t|]\\
&\leq M \E\left[\int_0^T \left(|\hat{U}_s|+|\hat{V}_s|+\E[|\hat{U}_s|]+\E[|\hat{V}_s|]\right)d s\right ].
\end{align*}
Similar estimates hold for $|s^1_t-s^2_t|$ and $|r^1(t,x)-r^2(t,x)|$. The above analysis yields that 
\begin{align}
 \sup_{t\in[0,T]}|\hat{K}_t|\leq M \E\left[\int_0^T \left(|\hat{U}_s|+|\hat{V}_s|+\E[|\hat{U}_s|]+\E[|\hat{V}_s|]\right)d s \right]. 
\end{align}
Applying H{\"o}lder inequality, we obtain 
\begin{align}\label{differK}
\sup_{t\in[0,T]}|\hat{K}_t|^p\leq M\E\left[\left(\int_0^T \left(|\hat{U}_s|+|\hat{V}_s|+\E[|\hat{U}_s|]+\E[|\hat{V}_s|]\right)ds\right)^p \right].
\end{align}
Combining Eqs. \eqref{differY} and \eqref{differK} yields that 
\begin{align*}
\E\left[\sup_{t\in[0,T]}|\hat{Y}_t|^p\right]\leq M\E\left[\left(\int_0^T \left(|\hat{U}_s|+|\hat{V}_s|+\E[|\hat{U}_s|]+\E[|\hat{V}_s|]\right)d s \right)^p\right].
\end{align*}

Note that we have
\begin{align*}
\int_0^t \hat{Z}_s dB_s=\int_0^t\hat{f}_s ds+\hat{Y}_t-\hat{Y}_0+\hat{K}_t.
\end{align*}
Simple calculation yields that
\begin{align*}
\E\left[\sup_{t\in[0,T]}\left(\int_0^t \hat{Z}_sd B_s\right)^p\right]&\leq M\E\left[\sup_{t\in[0,T]}\left(\int_0^t\hat{f}_s ds\right)^p+\sup_{t\in[0,T]}\left(\hat{Y}_t\right)^p+\sup_{t\in[0,T]}\left(\hat{K}_t\right)^p\right]\\
&\leq M\E\left[\left(\int_0^T \left(|\hat{U}_s|+|\hat{V}_s|+\E[|\hat{U}_s|]+\E[|\hat{V}_s|]\right)ds\right)^p\right],
\end{align*}
which, together with the B-D-G inequality yields that 
\begin{align*}
   \E\left[\left(\int_0^T |\hat{Z}_s|^2 d s\right)^{\frac{p}{2}}\right]&\leq M\E\left[\left(\int_0^T \left(|\hat{U}_s|+|\hat{V}_s |+\E[|\hat{U}_s|]+\E[|\hat{V}_s|]\right)ds\right)^p\right].
\end{align*}
By the H\"{o}lder inequality and the Fubini theorem, we finally deduce that 
\begin{align*}
\E\left[\sup_{s\in[0,T]}|\hat{Y}_s|^p+\left(\int_0^T |\hat{Z}_s|^2 d s\right)^{\frac{p}{2}}\right]
\leq& M\E\left[\left(\int_0^T \left(|\hat{U}_s|+|\hat{V}_s|+\E[|\hat{U}_s|]+\E[|\hat{V}_s|]\right)ds\right)^p\right]\\
\leq& M\E\left[\left(\int_0^T |\hat{U}_s|+|\hat{V}_s|ds\right)^p\right]\\
\leq& M\max\left(T^{\frac{p}{2}},T^p\right)\E\left[\sup_{t\in[0,T]}|\hat{U}_t|^p+\left(\int_0^T|\hat{V}_s|^2ds\right)^{\frac{p}{2}}\right].
\end{align*}
Let $T$ be sufficiently small such that $M\max\left(T^{\frac{p}{2}},T^p\right)<1$, we have constructed a contraction mapping. Therefore, when $T$ small enough, the doubly mean reflected $\text{MFBSDE}$ \eqref{nonlinearyz} has a unique solution.

For the general $T$, similar as a standard $\text{BSDE}$ approach, we divide the time interval $[0,T]$ into finitely many small time intervals. On each small time interval, the previous analysis ensures the existence of a local solution. By stitching the local solutions, we get a global solution on the whole time interval. Uniqueness on the whole interval is a direct consequence of the uniqueness on each small interval. The proof is complete. 
\end{proof}

\begin{remark}
    Another frequently used method for constructing the solution to the reflected problem is approximation via penalization (see \cite{CK, EKPPQ} for the classical reflected $\text{BSDEs}$, \cite{BEH,L} for the mean reflected $\text{BSDEs}$ and \cite{Li'} for the  reflected $\text{MFBSDEs}$). If the loss functions $L,R$ are linear, the penalization method is valid for the $\text{MFBSDE}$ with double mean reflections. The detailed proof is postponed to the Appendix.
\end{remark}

\section{Doubly mean reflected \text{MFBSDEs} with quadratic generator}	
	In this section, we investigate the doubly mean reflected \text{MFBSDE} \eqref{nonlinearyz} with quadratic generator $f$. Since $f$ is not Lipschitz continuous, the method in the proof of Theorem \ref{th-1} fails. We will consider the bounded terminal condition and unbounded terminal condition separately. For the first case, we use the $\text{BMO}$ martingale theory and a fixed-point argument. For the second case, we need to apply the $\theta$-method and introduce either convexity or concavity on the generator. Besides, in both cases, the generator does not depend on the distribution of the $Z$-term. 
   \subsection{Bounded terminal value}
   To investigate the existence and uniqueness of solution to the doubly mean reflected \text{MFBSDE} \eqref{nonlinearyz} with quadratic generator $f$ and bounded terminal condition $\xi$,  we need introduce the following assumption.
   \begin{assumption}\label{ass-3}
  The process $\left(f\left(t, 0, \delta_0, 0\right)\right)$ is uniformly bounded and there exist two positive constants $\beta$ and $\gamma$ such that for any $t \in[0, T], y_1, y_2 \in \mathbb{R}, v_1, v_2 \in \mathcal{P}_1(\mathbb{R})$, and $z_1, z_2 \in \mathbb{R}$.
	$$
    \begin{aligned}
		& \left|f\left(t, y_1, v_1, z_1\right)-f\left(t, y_2, v_2, z_2\right)\right| \\
		& \left.\leq \lambda\left(\left|y_1-y_2\right|+W_1\left(v_1, v_2\right)\right)+\gamma\left(1+\left|z_1\right|+\left|z_2\right|\right)\right)\left|z_1-z_2\right| .
	\end{aligned}
	$$
 \end{assumption}

	\begin{theorem}\label{thmbddterminal}
	 Suppose that  $\xi \in \mathcal{L}_T^{\infty}$ satisfies $\mathbf{E}[L(T,\xi)]\leq 0 \leq \mathbf{E}[R(T,\xi)]$ and Assumptions \ref{ass2}, \ref{ass-3} hold. Then the quadratic MFBSDE \eqref{nonlinearyz} with double mean reflections admits a unique solution $(Y, Z, K) \in \mathcal{S}^{\infty} \times \text{BMO}\times BV[0,T]$.
	\end{theorem}
\begin{proof}
The proof will be divided into the following three parts.

{\bf Step 1.} We claim that given $U\in \mathcal{S}^\infty$, the following $\text{MFBSDE}$ with double mean reflections:
\begin{align}\label{MFBSDEDMR-1}
\begin{cases}
	Y^U_t=\xi+\int_t^T f\left(s,U_{s}, \P_{U_{s}}, Z^U_{s}\right)ds-\int_t^T Z^U_s dB_s+\left(K^U_T-K^U_t\right),\ t\in[0,T],\\
	\E\left[L(t,Y^U_t)\right]\leq 0 \leq \E\left[R(t,Y^U_t)\right], \ K^U_t=K_t^{U,R}-K_t^{U,L},\ K_t^{U,R},K_t^{U,L}\in I[0,T],  \\
\int_{0}^{T}\E\left[R(t,Y^U_t)\right]dK_t^{U,R}=\int_{0}^{T}\E\left[L(t,Y^U_t)\right]dK_t^{U,L}=0.\\
\end{cases}
\end{align}
admits a unique solution $(Y^U,Z^U,K^U)\in \mathcal{S}^{\infty} \times \text{BMO}\times BV[0,T]$. 
For the uniqueness, suppose that $(Y^i,Z^i,K^i)$ are the solutions to \eqref{MFBSDEDMR-1}, $i=1,2$. Since $(Y^i-(K^i_T-K^i),Z^i)$ can be seen as the solution to the $\text{BSDE}$ with terminal value $\xi$ and generator $\{f(s,U_s,\P_{U_s},z)\}_{s\in[0,T]}$, by Theorem 7.3.3 in \cite{ZJ}, we have $Y^1-(K^1_T-K^1)\equiv Y^2-(K^2_T-K^2)$ and $Z^1\equiv Z^2(=:z)$. It remains to show $K^1\equiv K^2$. Set $C_s=f(s,U_s,\P_{U_s}, Z_s)$, the proof of $K^1\equiv K^2$ is the same as the one of Proposition 4.2 in \cite{L}. So we omit it. 

For the existence, by Theorem 7.3.3 in \cite{ZJ}, the following $\text{BSDE}$
\begin{align*}
    y^U_t=\xi+\int_t^T f(s,U_s,\P_{U_s},z^U_s)ds-\int_t^T z^U_s dB_s
\end{align*}
admits a unique solution $(y^U,z^U)\in \mathcal{S}^\infty\times \text{BMO}$. Set 
\begin{equation}\label{salr}\begin{split}
&s^U_t=\E\left[\int_0^t f(s,U_s,\P_{U_s},z^U_s) ds\right]=\E\left[y^U_0-y^U_t\right], \ a=\E\left[\xi\right],\\
&l^U(t,x):=\E\left[L(t,y^U_t-\E\left[y^U_t\right]+x)\right], \ r^U(t,x):=\E\left[R(t,y^U_t-\E\left[y^U_t\right]+x)\right].
\end{split}\end{equation}
It is easy to check that $l^U,r^U$ satisfy Assumption \ref{asslr} and 
\begin{align*}
l^U(T,a)=\E\left[L(T,\xi)\right]\leq 0\leq \E\left[R(T,\xi)\right]=r^U(T,a).
\end{align*}
Therefore, by Theorem \ref{BSP}, the backward Skorokhod problem $\mathbb{BSP}_{l^U}^{r^U}(s^U,a)$ admits a unique solution $(x^U,K^U)$ and for any $t\in[0,T]$, $K^U_t=\bar{K}^U_T-\bar{K}^U_{T-t}$, where $\bar{K}^U$ is the second component of the solution to the Skorokhod problem $\mathbb{SP}_{\bar{l}^U}^{\bar{r}^U}(\bar{s}^U)$. Here, for any $t\in[0,T]$, 
\begin{equation}\label{barslr}\bar{s}^U_{t}=a+s^U_{T}-s^U_{T-t}=\mathbf{E}[y^U_{T-t}],\ \bar{l}^U(t,x)=l^U (T-t,x), \ \bar{r}^U(t,x)=r^U (T-t,x).
\end{equation}
Now, we set 
\begin{align}\label{re-1}
Y^U_t=y^U_t+K^U_T-K^U_t=y^U_t+\bar{K}^U_{T-t}, \ Z^U_t=z^U_t, \ t\in[0,T].
\end{align}
Then, $(Y^U,Z^U,K^U)$ is the solution to \eqref{MFBSDEDMR-1}, whose proof is the same as the one of  Proposition 4.2 in \cite{L}.

{\bf Step 2.} By Step 1, we define a map $\Gamma:\mathcal{S}^\infty\rightarrow \mathcal{S}^\infty$ as follows
\begin{equation}
\Gamma(U):=Y^U,
\end{equation}
where $Y^U$ is the first component of the solution to \eqref{MFBSDEDMR-1}. We claim that $\Gamma$ is a contraction mapping when the terminal time $T$ is sufficiently small.

Given $U^{i}\in\mathcal{S}^\infty$, $i=1,2$, let $(y^i,z^i)\in \mathcal{S}^\infty\times \text{BMO}$ be the solution to the following BSDE:
$$y_t^i=\xi+\int_t^T f(s,U_s^i,\P_{U_s^i}, z^i_s)ds-\int_t^T z^i_s dW_s.$$
Define $\hat{y}_{t}:=y^1_t-y^2_t,\hat{z}_{t}:=z^1_t-z^2_t,\hat{U}_{t}:=U^1_t-U^2_t$. For each $t \in[0, T]$, we denote
\begin{equation}\label{re-7}
   \beta_t=\frac{f\left(t, U_t^1,\P_{U_t^1},z_t^1\right)-f\left(t, U_t^1,\P_{U_t^1}, z_t^2\right)}{|\hat{z}_t|^2}\hat{z}_t \mathbf{1}_{\left\{|\hat{z}_t| \neq 0\right\}} . 
\end{equation}
It is easy to check that the pair of processes $(\hat{y}, \hat{z})$ satisfies the following $\text{BSDE}$:
$$
\begin{aligned}
	\hat{y}_t=  \int_t^T\left(\beta_s \hat{z}_s+f\left(s,U_s^1,\P_{U_s^1}, z_s^2\right)-f\left(s,U_s^2,\P_{U_s^2}, z_s^2\right)\right) d s 
	 -\int_t^T \hat{z}_s d B_s .
\end{aligned}
$$
Solving the above $\text{BSDE}$ yields that $$
\hat{y}_t=\mathbf{E}_t\left[\mathscr{E}(\beta \cdot B)_t^T\left(\int_t^T\left(f\left(s,U_s^1,\P_{U_s^1}, z_s^2\right)-f\left(s,U_s^2,\P_{U_s^2}, z_s^2\right)\right) d s\right)\right],
$$ 
where for any $0\leq u\leq v\leq T$, 
\begin{align*}
    \mathscr{E}\left(\beta \cdot B\right)_u^v=\exp\left\{\left(\int_u^v \beta_s dB_s-\frac{1}{2}\int_u^v \left|\beta_s\right|^2ds\right)\right\}.
\end{align*}
 Applying Assumption $\ref{ass-3}$, we know $\beta\in \text{BMO}$, which implies that $\left\{\mathscr{E}\left(\beta \cdot B\right)_0^t\right\}_{t\in[0,T]}$ is a martingale. Consequently, for any $t \in[0, T]$, we have 
\begin{equation}\label{hatyt}\begin{split}
\left|\hat{y}_t\right|&\leq\mathbf{E}_t\left[\left|\mathscr{E}(\beta \cdot B)_t^T\right| \left(\int^T_0 \left|f\left(s,U_s^1,\P_{U_s^1}, z_s^2\right)-f\left(s,U_s^2,\P_{U_s^2}, z_s^2\right)\right|ds\right)\right]\\
 &\leq\mathbf{E}_t\left[|\mathscr{E}(\beta \cdot B)_t^T|\lambda \left(\int^T_0 \left(|\hat{U}_{s}|+W_1(P_{U_s^1},P_{U_s^2})\right)ds\right)\right]\\
 &\leq\left(\lambda T\|\hat{U}\|_{ \mathcal{S}^{\infty}}+\lambda T \sup _{s \in[0,T]} \mathbf{E}[|\hat{U}_s|]\right)\mathbf{E}_t\left[\mathscr{E}(\beta \cdot B)_t^T\right]\\
 &\leq 2\lambda T\|\hat{U}\|_{ \mathcal{S}^{\infty}} . 
\end{split}
\end{equation}

For notional simplification, set $s^i=s^{U^i}$, $l^i=l^{U^i}$, $r^i=r^{U^i}$ and $\bar{s}^i$, $\bar{l}^i$, $\bar{r}^i$ are defined similar as \eqref{barslr}, $i=1,2$. Let $\bar{K}^i$ be the solution to the Skorokhod problem $\mathbb{SP}_{\bar{l}^i}^{\bar{r}^i}(\bar{s}^i)$, $i=1,2$. In view of Proposition \ref{continuity}, we know
$$
	\sup_{t \in[0, T]}\left|\bar{K}^1_{t}-\bar{K}^2_{t}\right|\leq \frac{C}{c}\sup_{t\in[0,T]}\left|\bar{s}^1_t-\bar{s}^2_t\right|+\frac{1}{c}(\bar{L}_T \vee \bar{R}_T),	
$$
where $\bar{L}_T=\sup_{(t,x)\in[0,T]\times \mathbb{R}}|\bar{l}^1(t,x)-\bar{l}^2(t,x)|$, $
\bar{R}_T=\sup_{(t,x)\in[0,T]\times \mathbb{R}}|\bar{r}^1(t,x)-\bar{r}^2(t,x)|$. Recalling that  $\bar{s}^i_t=\E[y^i_{T-t}]$, we have 
\begin{align*}
	\sup_{t\in[0,T]}|\bar{s}^1_t-\bar{s}^2_t|
    \leq \sup_{t\in[0,T]}\mathbf{E}\left[\left|\hat{y}_t\right|\right].
	\end{align*}
Applying Assumption \ref{ass2}, we obtain that 
$$
\begin{aligned}
	\bar{L}_T
	&= \sup_{(t,x)\in[0,T]\times \mathbb{R}} \left|\mathbf{E}\left[L(T-t,y^1_{T-t}-\mathbf{E}\left[y^1_{T-t}\right]+x)\right]-\E\left[L(T-t,y^2_{T-t}-\mathbf{E}\left[y^2_{T-t}\right]+x)\right]\right|\\
	&\leq \sup_{t\in[0,T]}\mathbf{E}\left[C\left|\hat{y}_{T-t}-\mathbf{E}\left[\hat{y}_{T-t}\right]\right|\right]
	\leq 2C\sup_{t \in[0, T]}\mathbf{E}\left[\left|\hat{y}_t\right|\right].
\end{aligned}
$$
Similarly, we have $\bar{R}_T\leq 2C\sup_{t \in[0, T]}\mathbf{E}[|\hat{y}_t|]$. The above analysis implies that 
\begin{equation}\label{hatbarK}
	\sup_{t \in[0, T]}\left|\bar{K}^1_{t}-\bar{K}^2_{t}\right|
	\leq \frac{3C}{c}\sup_{t \in[0, T]}\mathbf{E}\left[|\hat{y}_t|\right].
\end{equation}
Then, according to \eqref{re-1}, \eqref{hatyt} and  \eqref{hatbarK}, we can derive that
$$
\begin{aligned}
\left\|\Gamma\left(U^1\right)-\Gamma\left(U^2\right)\right\|_{\mathcal{S}^{\infty}}&=\left\| \hat{y}+\bar{K}^1-\bar{K}^2\right\|_{\mathcal{S}^{\infty}}\\
&\leq \left\| \hat{y}\right\|_{\mathcal{S}^{\infty}}+\sup_{t\in[0,T]}\left|\bar{K}^1_t-\bar{K}^2_t\right| \\
&\leq M(c,C,\lambda)T\|\hat{U}\|_{\mathcal{S}^{\infty}},\\ 
\end{aligned}
$$
where $M(c,C,\lambda)=2(1+\frac{3C}{c})\lambda$.

Therefore, if $T$ is sufficiently small such that  $M(c,C,\lambda)T<1$, $\Gamma$ is a contraction mapping. In this case, if the mapping $\Gamma$ has a unique fixed point $Y \in \mathcal{S}^\infty$,  the mean reflected $\text{BSDE}$ \eqref{MFBSDEDMR-1} with generator $f(\cdot, Y_{\cdot}, \P_{Y_{\cdot}}, \cdot)$   has a unique solution $(\widetilde{Y},Z,K)\in\mathcal{S}^{\infty} \times \text{BMO}\times BV[0,T]$. And we have $\widetilde{Y}=\Gamma(Y)=Y$, so $(Y,Z,K)$ is the unique solution to the mean field $\text{BSDE}$ \eqref{nonlinearyz}. For the general $T$, similar to the analysis of Theorem \ref{th-1}, we can obtain existence and uniqueness on the whole time interval.  The proof is complete. 
\end{proof}

\subsection{Unbounded terminal value} 
In this subsection, we consider the doubly mean reflected MFBSDE \eqref{nonlinearyz} with quadratic generator $f$ and unbounded terminal value $\xi$. The well-posedness can be established under the following conditions on the generator $f$.
\begin{assumption}\label{ass-4}
There exists a positive progressively measurable process $\left(\alpha_t\right)_{0 \leq t \leq T}$ with $\int_0^T \alpha_t d t \in \mathbb{L}_T$ and two positive constants $\beta$ and $\gamma$ such that
\begin{itemize}
\item[1.] $\forall(t, y, v, z) \in[0, T] \times \mathbb{R} \times \mathcal{P}_1(\mathbb{R}) \times \mathbb{R}$,
$$
|f(t, y, v, z)| \leq \alpha_t+\lambda \left(|y|+W_1\left(v, \delta_0\right)\right)+\frac{\gamma}{2}|z|^2,
$$
\item[2.] $\forall t \in[0, T], y_1, y_2 \in \mathbb{R}, v_1, v_2 \in \mathcal{P}_1(\mathbb{R}), z \in \mathbb{R}$,
$$
\left|f\left(t, y_1, v_1, z\right)-f\left(t, y_2, v_2, z\right)\right| \leq \lambda \left(\left|y_1-y_2\right|+W_1\left(v_1, v_2\right)\right),
$$
\item[3.] $\forall(t, y, v) \in[0, T] \times \mathbb{R} \times \mathcal{P}_1(\mathbb{R}),  z \mapsto f(t, y, v, z)$ is convex or concave.
\end{itemize}
\end{assumption}

\begin{theorem}\label{th-3}
Suppose that  $\xi \in \mathbb{L}_T$ satisfies $\E[L(T,\xi)]\leq 0 \leq \E[R(T,\xi)]$ and Assumptions $\ref{ass2}$, $\ref{ass-4}$ hold. Then the quadratic  MFBSDE  with double mean reflections \eqref{nonlinearyz} admits at most one solution $(Y, Z, K) \in \mathbb{S} \times \mathcal{H}^p \times BV[0,T]$.
\end{theorem}

Applying Corollary 6 in \cite{BH08}, Step 1 in the proof of Theorem \ref{thmbddterminal} is valid when $f$ satisfies Assumption \ref{ass-4}. More precisely, the representation \eqref{re-1} still holds. However, under Assumption $\ref{ass-4}$, the process $\eqref{re-7}$ may be unbounded in the $\text{BMO}$ space. Therefore,  $\{\mathscr{E}(\beta \cdot B)_0^t\}_{t\in[0,T]}$ is not a martingale. Then the fixed point argument, especially \eqref{hatyt}, fails to work. Motivated by \cite{HMW}, we use the $\theta$-method to prove the uniqueness and existence. Before proving Theorem $\ref{th-3}$, we first introduce the following technical lemma, which will be frequently used in the sequel.
\begin{lemma}\label{le-2}
 For any given finite time interval $[u,v]\subset [0,T]$ and any positive constants $C_i$, $i=0,1,\cdots, 6$,  suppose that $\xi \in \mathbb{L}_v$, $A, a, A'\in \mathbb{S}[u,v]$ are nonnegative and  satisfy the following conditions: 
 \begin{itemize}
 \item[(i)] $A_t\leq a_t+C_0\mathbf{E}\left[\sup_{t\in[u,v]}a_t\right]+C_1$,
 \item[(ii)]
$$\mathbf{E}\left[\exp\left\{C_2\sup_{t\in[u,v]}a_t\right\}\right]\leq C_3\mathbf{E}\left[\exp\left\{ C_2 C_4\xi \right\}\right]^{C_5}\mathbf{E}\left[\exp\left\{C_2 C_6\left(v-u\right)\sup_{t\in[u,v]}A'_t\right\}\right].$$
\end{itemize}
Then, for any $C_7>0$, we have 
 $$
 \begin{aligned}
  \mathbf{E}\left[\exp\left\{C_7\sup_{t\in[u,v]}A_t\right\}\right]&\leq C_3\exp\left\{C_1C_7\right\}\mathbf{E}\left[\exp\left\{ 2\xi\left(1+C_0\right)C_4C_7\right\}\right]^{C_5}\\
 &\times \mathbf{E}\left[\exp\left\{2\left(1+C_0\right)C_6C_7\left(v-u\right)\sup_{t\in[u,v]}A'_t\right\}\right].\\   
 \end{aligned}
$$  
\end{lemma}
\begin{proof}
$$
\begin{aligned}
 &\mathbf{E}\left[\exp\left\{C_7\sup_{t\in[u,v]}A_t\right\}\right]\\
 \leq & \mathbf{E}\left[\exp\left\{C_7\sup_{t\in[u,v]}a_t+C_0C_7\mathbf{E}\left[\sup_{t\in[u,v]}a_t\right]+C_1C_7\right\}\right]\\
 \leq&\exp\left\{C_1C_7\right\}\mathbf{E}\left[\exp\left\{C_7\sup_{t\in[u,v]}a_t\right\}\right]\mathbf{E}\left[\exp\left\{C_0C_7\sup_{t\in[u,v]}a_t\right\}\right]\\
 \leq&\exp\left\{C_1C_7\right\}\mathbf{E}\left[\exp\left\{(C_7+C_0C_7)\sup_{t\in[u,v]}a_t\right\}\right]\mathbf{E}\left[\exp\left\{(C_7+C_0C_7)\sup_{t\in[u,v]}a_t\right\}\right]\\
 \leq &\exp\left\{C_1C_7\right\}\mathbf{E}\left[\exp\left\{2(C_7+C_0C_7)\sup_{t\in[u,v]}a_t\right\}\right]\\
 \leq& C_3 \exp\left\{C_1C_7\right\}\mathbf{E}\left[\exp\left\{2\xi (C_7+C_0C_7) C_4\right\}\right]^{C_5}\mathbf{E}\left[\exp\left\{2(C_7+C_0C_7)C_6\left(v-u\right)\sup_{t\in[u,v]}A'_t\right\}\right].\\
\end{aligned}
$$
where we have used condition (i) in the first inequality,  Jensen's inequality in the second inequality,  H\"{o}lder's inequality in the fourth inequality and condition (ii) in the last inequality.
\end{proof}

\begin{proof}[Proof of Theorem \ref{th-3}] Without loss of generality, we assume that $f(t, y, v, \cdot)$ is concave. For $i=1,2$, let  $\left(Y^i, Z^i, K^i\right)$ be the solutions to the quadratic mean reflected $\text{MFBSDE}$  \eqref{nonlinearyz}. For each $\theta \in(0,1)$, we define
\begin{equation}\label{deltay}
\delta_\theta Y:=\frac{\theta Y^1-Y^2}{1-\theta}, \quad \delta_\theta \tilde{Y}:=\frac{\theta Y^2-Y^1}{1-\theta} \quad \text { and } \quad \delta_\theta \bar{Y}:=\left|\delta_\theta Y\right|+|\delta_\theta \tilde{Y}|.
\end{equation}
 For any $\gamma>0$, it is easy to check that 
\begin{equation}\label{esupY}\begin{split}
	\mathbf{E} \left[\sup _{t \in[0, T]}\left|Y_t^1-Y_t^2\right|\right]&=\mathbf{E} \left[\sup _{t \in[0, T]}\left|(1-\theta)(\delta_{\theta}Y_{t}+Y^1_t)\right|\right] \\
  &\leq (1-\theta)\mathbf{E}\left[\sup_{t \in[0, T]}\left|\delta_{\theta}Y_{t}\right|+\sup _{t \in[0, T]}\left|Y^1_t\right| \right]\\
  &\leq (1-\theta)\mathbf{E}\left[\sup_{t \in[0, T]}\delta_\theta \bar{Y}_t+\sup _{t \in[0, T]}\left|Y^1_t\right| \right]\\
  & \leq(1-\theta)\left(\frac{1}{\gamma} \sup _{\theta \in(0,1)} \mathbf{E}\left[\exp \left\{\gamma \sup _{t \in[0, T]} \delta_\theta \bar{Y}_t\right\}\right]+\mathbf{E}\left[\sup _{t \in[0, T]}\left|Y_t^1\right|\right]\right).
  \end{split}
\end{equation}
Providing that 
\begin{equation}\label{bddexpdetltay}
    \sup _{\theta \in(0,1)} \mathbf{E}\left[\exp \left\{\gamma \sup _{t \in[0, T]} \delta_\theta \bar{Y}_t\right\}\right]<\infty,
\end{equation}
letting $\theta \to 1$ in \eqref{esupY}, we obtain that $Y^1=Y^2$. Then, applying Itô's formula to $|Y^1-Y^2|^2$, we have $Z^1=Z^2,K^1=K^2$ on $[0,T]$. In the following, we prove that \eqref{bddexpdetltay} holds when $T$ is sufficiently small, which will be divided into the following three parts.

{\bf Step 1.} We claim that  for any $\gamma>0$ and  $p> 1$, 
\begin{equation}\label{re-5}
\mathbf{E}\left[\exp \left\{p \gamma \sup _{t \in[0, T]} \delta_\theta \bar{y}_t\right\}\right]   \leq 4\mathbf{E}\left[\exp\left\{16p\gamma\left(|\xi|+\tilde{\chi}\right)\right\}\right]^{\frac{1}{2}}\mathbf{E}\left[\exp\left\{8p\gamma \lambda T \sup _{t \in[0, T]} \delta_\theta \bar{Y}_t\right\}\right],
\end{equation}
where $$\tilde{\chi}=  \int_0^T \alpha_s d s+2 \lambda T \sup_{t\in[0,T]} \mathbf{E}\left[\sum_{i=1}^2\left|Y^i_t\right|\right] +2 \lambda T\sum_{i=1}^2\sup _{t \in[0, T]}\left|Y_t^i\right| +\sum_{i=1}^2\sup _{t \in[0, T]}\left|y_t^i\right|$$ 
and  $\left(y^i, z^i\right) \in \mathbb{S} \times \mathcal{H}^p$ is the solution to the following quadratic $\text{BSDE}$:
$$
y_t^i=\xi+\int_t^T f\left(s, Y_s^i, \mathbf{P}_{Y_s^i}, z_s^i\right) d s-\int_t^T z_s^i d B_s .
$$
Here and in the sequel, $\delta_\theta \ell$, $\delta_\theta \tilde{\ell}$ and $\delta_\theta \bar{\ell}$ are defined simiarly as \eqref{deltay} for $\ell=y,z$. 

Indeed, simple calculation yields that 
$$
\begin{aligned}
 &\mathbf{E}\left[\exp \left\{p \gamma \sup _{t \in[0, T]} \delta_\theta \bar{y}_t\right\}\right]\\
 \leq& 4 \mathbf{E}\left[\exp \left\{4 p \gamma\left(|\xi|+\tilde{\chi}+\lambda T\left(\sup _{t \in[0, T]} \delta_\theta \bar{Y}_t+\sup _{t \in[0, T]} \mathbf{E}\left[\delta_\theta \bar{Y}_t\right]\right)\right)\right\}\right]\\
 \leq &4\mathbf{E}\left[\exp\left\{4p\gamma\left(|\xi|+\tilde{\chi}+\lambda T\sup _{t \in[0, T]} \delta_\theta \bar{Y}_t\right)\right\}\right]\mathbf{E}\left[\exp\left\{4p\gamma \lambda T \sup _{t \in[0, T]} \delta_\theta \bar{Y}_t\right\}\right]\\
 \leq& 4\mathbf{E}\left[\exp\left\{8p\gamma\left(|\xi|+\tilde{\chi}\right)\right\}\right]^{\frac{1}{2}}\mathbf{E}\left[\exp\left\{8p\gamma \lambda T \sup _{t \in[0, T]} \delta_\theta \bar{Y}_t\right\}\right],\\
\end{aligned}
$$
where we have used Eq. (17) in \cite{HMW} in the first inequality, Jensen's inequality in the second inequality and H\"{o}lder's inequality in the third inequality.

{\bf Step 2.} We show that for any $t\in[0,T]$,
\begin{equation}\label{deltathetaY}
    \begin{split}
    &|\delta_\theta Y_t| \leq B_1+|\delta_\theta y_t|+\frac{3C}{c} \sup _{0 \leq t \leq T} \mathbf{E}[|\delta_\theta y_t|],\\
	&|\delta_\theta \tilde{Y}_t| \leq B_1+|\delta_\theta \tilde{y}_t|+\frac{3C}{c} \sup _{0 \leq t \leq T} \mathbf{E}[|\delta_\theta \tilde{y}_t|],
\end{split}
\end{equation}
where 
$B_1:=\sup _{t \in[0, T]}|\bar{K}^0_{t}|+\frac{6C}{c} \sup _{t \in[0, T]} \mathbf{E}\left[\sum_{i=1}^2|y_t^i|\right]+\frac{C}{c}\E\left[|\xi|\right]$ and $\bar{K}^0$ is the second component of the solution to the Skorokhod problem $\mathbb{SP}_{\bar{l}^0}^{\bar{r}^0}(\bar{s}^0)$. Here $$
\bar{l}^0(t,x)=l^0(T-t,x)=\mathbf{E}\left[L(T-t,x)\right], \bar{r}^0(t,x)=r^0(T-t,x)=\mathbf{E}\left[R(T-t,x)\right], \bar{s}^0=\mathbf{E}[\xi].$$
We only need to show the first inequality in \eqref{deltathetaY} since the second one can be proved similarly.

Recalling \eqref{salr}, set $s^i=s^{Y^i}$, $l^i=l^{Y^i}$, $r^i=r^{Y^i}$, $i=1,2$ and $a=\E[\xi]$. Define $\bar{s}^i$, $\bar{l}^i$, $\bar{r}^i$, $i=1,2$, as in \eqref{barslr}. Let $\bar{K}^i$ be the second component of the solution to the Skorokhod problem $\mathbb{SP}_{\bar{l}^i}^{\bar{r}^i}(\bar{s}^i)$, $i=1,2$. Then, we have
$$
Y_t^i=y_t^i+\bar{K}^{i}_{T-t}, \quad \forall t \in[0, T].
$$
Consequently, we obtain that 
$$
\begin{aligned}
\left|\delta_{\theta}Y_{t}\right|&=\left|\frac{\theta\left(y_{t}^{1}+\bar{K}^{1}_{T-t}\right)-\left(y_{t}^{2}+\bar{K}^{2}_{T-t}\right)}{1-\theta}\right|\\
	&\leq \left|\delta_{\theta}y_{t}\right|+\frac{\theta}{1-\theta}\sup_{t\in[0,T]}\left|\bar{K}^{1}_{t}-\bar{K}^{2}_{t}\right|+\sup_{t\in[0,T]}\left|\bar{K}^{2}_{t}-\bar{K}^{0}_{t}\right|+\sup_{t\in[0,T]}\left|\bar{K}^{0}_{t}\right|.
\end{aligned}$$
By Proposition \ref{continuity}, we have
\begin{equation}\label{re-6}
 \sup_{t\in[0,T]}\left|\bar{K}^{2}_{t}-\bar{K}^{0}_{t}\right|\leq\frac{C}{c}\sup_{t\in[0,T]}\left|\bar{s}^{2}_{t}-\bar{s}^{0}_{t}\right|+\frac{1}{c}\left(\bar{L}^0_T\vee \bar{R}^0_T\right),
\end{equation}
where 
\begin{align*}
\bar{L}_T^0=\sup_{(t,x)\in[0,T]\times \mathbb{R}}|\bar{l}^2(t,x)-\bar{l}^0(t,x)|, \ 
\bar{R}_T^0=\sup_{(t,x)\in[0,T]\times \mathbb{R}}|\bar{r}^2(t,x)-\bar{r}^0(t,x)|.
\end{align*}
Recalling the definition of $\bar{s}^i$, $\bar{l}^i$, $\bar{r}^i$ for $i=0,2$, we derive that
\begin{align*}
\sup_{t\in[0,T]}\left|\bar{K}^{2}_{t}-\bar{K}^{0}_{t}\right| \leq \frac{3C}{c}\sup_{t\in[0,T]}\mathbf{E}\left[\left|y_t^2\right| \right]+\frac{C}{c}\mathbf{E}[\left|\xi\right|].
\end{align*}
By a similar analysis as  \eqref{hatbarK}, we obtain that
$$
\begin{aligned}
 \frac{\theta}{1-\theta}\sup_{t\in[0,T]}\left|\bar{K}^{1}_{t}-\bar{K}^{2}_{t}\right|&\leq \frac{\theta}{1-\theta}\frac{3C}{c}\sup_{t \in[0, T]}\mathbf{E}[|y_t^1-y_t^2|]\\
 &\leq \frac{3C}{c}\sup_{t\in[0,T]}\mathbf{E} \left[\left|\frac{\theta y_t^1-y_t^2+(1-\theta) y_t^2 }{1-\theta}\right|\right]\\
 &\leq \frac{3C}{c}\sup_{t\in[0,T]}\mathbf{E}\left[|\delta_{\theta}y_{t}|\right]+\frac{3C}{c}\sup_{t\in[0,T]}\mathbf{E}\left[|y_t^2|\right].\\   
\end{aligned}
$$
All the above analysis indicates that 
$$
\begin{aligned}
    |\delta_\theta Y_t| 
    &\leq|\delta_\theta y_t|+\frac{3C}{c} \sup _{0 \leq t \leq T} \mathbf{E}\left[|\delta_\theta y_t|\right]+\frac{6C}{c} \sup _{0 \leq t \leq T} \mathbf{E}[|y_t^2|]+\frac{C}{c}\mathbf{E}[|\xi|]+\sup_{t\in[0,T]}\left|\bar{K}^{0}_{t}\right|\\
    &\leq |\delta_\theta y_t|+\frac{3C}{c} \sup _{0 \leq t \leq T} \mathbf{E}[|\delta_\theta y_t|]+B_1.
\end{aligned}
$$

{\bf Step 3.} We are ready to prove \eqref{bddexpdetltay} holds when $T$ satisfies $(16+96\frac{C}{c})\lambda T<1$. First, by {\bf Step 2}, we obtain that
\begin{align*}
    \delta_\theta \bar{Y}_t&\leq 2B_1+\delta_\theta \bar{y}_t+\frac{6C}{c}\sup_{t\in[0,T]}\E[\delta_\theta \bar{y}_t]\\
    &\leq 2B_1+\delta_\theta \bar{y}_t+\frac{6C}{c}\E[\sup_{t\in[0,T]}\delta_\theta \bar{y}_t].
\end{align*}

Then, for any $\gamma>0$, $p > 1$ and $\theta \in (0,1)$, recalling \eqref{re-5}, by Lemma $\ref{le-2}$, we can obtain  
$$
\begin{aligned}
	\mathbf{E}\left[\exp \left\{p \gamma \sup _{t \in[0, T]} \delta_\theta \bar{Y}_t\right\}\right]
    &\leq \widetilde{C}_1\mathbf{E}\left[\exp\left\{\left(16+96\frac{C}{c}\right)p\gamma \lambda T\sup _{t \in[0, T]} \delta_\theta \bar{Y}_t\right\}\right]\\
	&\leq \widetilde{C}_1\mathbf{E}\left[\exp \left\{p \gamma \sup _{t \in[0, T]} \delta_\theta \bar{Y}_t\right\}\right]^{(16+96 \frac{C}{c}) \lambda T},
\end{aligned}
$$
where 
\begin{align*}
    \widetilde{C}_1=4\exp\left\{2p\gamma B_1\right\}\mathbf{E}\left[\exp\left\{\left(16+96\frac{C}{c}\right)p\gamma\left(|\xi|+\tilde{\chi}\right)\right\}\right]^{\frac{1}{2}}<\infty.
\end{align*}

Finally, we obtain that 
$$
\begin{aligned}
\mathbf{E}\left[\exp \left\{p \gamma \sup _{t \in[0, T]} \delta_\theta \bar{Y}_t\right\}\right]\leq \widetilde{C}_1^{\frac{1}{1-(16+96\frac{C}{c}\lambda T)}}<\infty.
\end{aligned}
$$
 Therefore \eqref{bddexpdetltay} holds true. Then, we have established the uniqueness result when $T$ is sufficiently small. For the general $T$, uniqueness on the whole interval is a direct consequence of the uniqueness on each small interval.
 \end{proof}
 
 Next, we aim to prove the existence of the solution to Eq. \eqref{nonlinearyz} whose generator is independent of the distribution of $Z$. We only need to consider the case on a sufficiently small time interval, since  for the general $T$, the construction can be obtained by  similar analysis as in the proof of Theorem \ref{th-1}. 
 \begin{assumption}\label{assl'r'}
     For any $(\omega,t,x)\in\Omega\times [0,T]\times\mathbb{R}$, we have 
     \begin{align*}
         L(\omega,t,x){\color{red}\leq} L'(t,x), \ R(\omega,t,x)\geq R'(t,x),
     \end{align*}
     where 
     \begin{align*}
         L'(t,x)=b_t x-p_t, \ R'(t,x)=b_tx-q_t.
     \end{align*}
     Here, $b,p,q$ are continuous functions such that
     \begin{itemize}
         \item $b$ is strictly positive;
         \item $\inf_{t\in[0,T]}(p_t-q_t)>0$;
         \item the functions $\frac{p}{b},\frac{q}{b}$ have bounded variation on $[0,T]$.
     \end{itemize}
 \end{assumption}
 \begin{theorem}\label{4-6}
Suppose that Assumptions \ref{ass2}, \ref{ass-4} and \ref{assl'r'} hold.  For any fixed $t_0\in[0,T)$, let $h\in(0,T-t_0]$ be such that $(32+192\frac{C}{c})\lambda h<1$  and let $\eta\in \mathbb{S}(\mathcal{F}_{t_0+h})$ be such that 
\begin{align*}
    \E[L(t_0+h,\eta)]\leq 0\leq \E[R(t_0+h,\eta)].
\end{align*}
Then the following MFBSDE with double mean reflections has
 unique solution $(Y,Z,K) \in \mathbb{S} [t_0,t_0+h]\times \mathcal{H}^p [t_0,t_0+h] \times BV[t_0,t_0+h]$ 
\begin{equation}\label{Yt0}
 \begin{cases}
	Y_t=\eta+\int_t^{t_{0}+h} f\left(s,Y_s, \mathbf{P}_{Y_s}, Z_s\right) d s-\int_t^{t_{0}+h} Z_s d B_s +K_{t_{0}+h}-K_t, \quad t_{0} \leq t \leq t_{0}+h, \\
	\E[L(t,Y_t)]\leq 0 \leq E[R(t,Y_t)],\ K_t=K_t^R-K_t^L,\ K_t^R,K_t^L\in I[t_{0},t_{0}+h],  \\
	\int_{t_{0}}^{t_{0}+h}\E[R(t,Y_t)]dK_t^R=\int_{t_{0}}^{t_{0}+h}\E[L(t,Y_t)]dK_t^L=0.
\end{cases}   
\end{equation}   
 \end{theorem}


\begin{proof}
Uniqueness can be obtained to use Theorem \ref{th-3}. For existence, set $Y^{(0)}=0$. Then, we define the stochastic process sequences $\left\{\left(Y^{(m)}, Z^{(m)}, K^{(m)}\right)\right\}_{m=1}^{\infty}$ recursively by solving  the following quadratic $\text{BSDE}$ with double mean reflections:
\begin{equation}\label{Ytm}
\begin{cases}
	Y_t^{(m)}=\eta+\int_t^{t_{0}+h} f\left(s, Y_s^{(m-1)}, \mathbf{P}_{Y_s^{(m-1)}}, Z^{(m)}_s\right) ds-\int_t^{t_{0}+h} Z_s^{(m)} dB_s +K_{t_{0}+h}^{(m)}-K_t^{(m)}, \\
    \E[L(t,Y^{(m)}_t)]\leq 0 \leq \E[R(t,Y^{(m)}_t)], \ t_{0} \leq t \leq t_{0}+h,\\ 
	K_t^{(m)}=(K_t^R)^{(m)}-(K_t^L)^{(m)},\ (K_t^R)^{(m)},(K_t^L)^{(m)}\in I[t_{0},t_{0}+h],  \\
\int_{t_{0}}^{t_{0}+h}\E[R(t,Y_t^{(m)})]d(K_t^R)^{(m)}=\int_{t_{0}}^{t_{0}+h}\E[L(t,Y_t^{(m)})]d(K_t^L)^{(m)}=0.\\
\end{cases}
\end{equation}
We show that $\left\{\left(Y^{(m)}, Z^{(m)}, K^{(m)}\right)\right\}_{m=1}^{\infty}$ converges to $(Y,Z,K)$, which is the solution to Eq. \eqref{Yt0}. The proof will be divided into the following parts.    

{\bf Step 1.} We give the representation for $Y^{(m)}$. The well-posedness of the above quadratic mean reflected $\text{MFBSDE}$ can be proved similarly to {\bf Step 1} in the proof of Theorem $\ref{thmbddterminal}$, where we just need to replace Theorem 4.2 in \cite{BDHPS} by Corollary 6 in \cite{BH08}. 
Let $(y^{(m)},z^{(m)})$ be the solution to the following quadratic $\text{BSDE}$ 
\begin{equation}\label{ytm}
y_t^{(m)}=\eta+\int_t^{t_0+h} f\left(s, Y_s^{(m-1)}, \mathbf{P}_{Y_s^{(m-1)}}, z^{(m)}_s\right) d s-\int_t^{t_0+h} z^{(m)}_s dB_s .    
\end{equation}
Since $K^{(m)}$ is a deterministic function, according to Corollary 6 in \cite{BH08}, we obtain that 
$$\left( Y^{(m)}_{\cdot}-(K_{t_0+h}^{(m)}-K_{\cdot}^{(m)}),\ Z_{\cdot}^{(m)} \right)\in\mathbb{S} [t_0,t_0+h]\times \mathcal{H}^p [t_0,t_0+h]$$
is again a  solution to the $\text{BSDE}$ \eqref{ytm}, which implies that \begin{equation}\label{Zmzm}\left(Y^{(m)}_t,Z^{(m)}_t\right)=\left(y^{(m)}_t+K^{(m)}_{t_0+h}-K^{(m)}_t,z^{(m)}_t\right),\quad \forall t\in [t_0,t_0+h].
\end{equation}

For any $(t,x)\in[t_0,t_0+h]\times \mathbb{R}$, we define 
\begin{equation}
\begin{split}
&s^{(m)}_t:=\E
\left[\int_{t_0}^t f\left(s,Y_s^{(m-1)}, \mathbf{P}_{Y_s^{(m-1)}},Z^{(m)}_s\right)ds\right]=\E\left[y^{(m)}_{t_0}-y^{(m)}_t\right], \ a^{(m)}:=\E[\eta],\\
&l^{(m)}(t,x):=\E\left[L\left(t,y^{(m)}_t-\E[y^{(m)}_t]+x\right)\right], \ r^{(m)}(t,x):=\E\left[R\left(t,y^{(m)}_t-\E[y^{(m)}_t]+x\right)\right].
\end{split}
\end{equation}
It is easy to check that $l^{(m)},r^{(m)}$ satisfy Assumption \ref{asslr} and 
\begin{align*}
l^{(m)}(t_0+h,a)=\E[L(t_0+h,y^{(m)}_{t_0+h})]\leq 0\leq \E[R(t_0+h,y^{(m)}_{t_0+h}))]=r^{(m)}(t_0+h,a).
\end{align*}
Therefore, by Theorem \ref{BSP}, the backward Skorokhod problem $\mathbb{BSP}_{l^{(m)}}^{r^{(m)}}(s^{(m)},a)$ admits a unique solution $(x^{(m)},K^{(m)})$. Furthermore, for any $t\in[t_0,t_0+h]$, we have $K^{(m)}_t=\bar{K}^{(m)}_h-\bar{K}^{(m)}_{t_0+h-t}$,  where $\bar{K}^{(m)}$ is the second component of the solution to the Skorokhod problem $\mathbb{SP}_{\bar{l}^{(m)}}^{\bar{r}^{(m)}}(\bar{s}^{(m)})$. Here, $\bar{s}^{(m)}$, $\bar{l}^{(m)}$ and $\bar{r}^{(m)}$ are given as follows 
\begin{equation}\label{barsbarlbarr}\begin{split}
&\bar{s}^{(m)}_{t}=a+s^{(m)}_{t_0+h}-s^{(m)}_{t_0+h-t}=\mathbf{E}[y^{(m)}_{t_0+h-t}],\\ 
&\bar{l}^{(m)}(t,x)=l^{(m)} (t_0+h-t,x), \ \bar{r}^{(m)}(t,x)=r^{(m)}(t_0+h-t,x).
\end{split}\end{equation}
 And for any $m\geq 1$, we  know    
\begin{equation}\label{re-2}
Y^{(m)}_t=y^{(m)}_t+K^{(m)}_{t_{0}+h}-K^{(m)}_t=y^{(m)}_t+\bar{K}^{(m)}_{t_{0}+h-t}, \quad \forall t\in[t_{0},t_{0}+h],   
\end{equation}
which is the desired representation for $Y^{(m)}$.

{\bf Step 2.}  We give some uniform estimates for $\left(Y^{(m)}, Z^{(m)}, K^{(m)}\right)\in\mathbb{S} [t_0,t_0+h]\times \mathcal{H}^p [t_0,t_0+h] \times BV[t_0,t_0+h]$. Show that for any $\gamma>0$ and $p >1$, we have
\begin{equation}\label{uniformbounded}
\sup _{m \geq 0} \mathbf{E}\left[\exp \left\{p \gamma \sup _{t \in[t_{0}, t_{0}+h]}|Y_t^{(m)}|\right\}+\left(\int_{t_{0}}^{t_{0}+h}|Z_t^{(m)}|^2 d t\right)^p+Var^{t_0+h}_{t_0}(K^{(m)})\right]<\infty 
\end{equation}
and
	\begin{equation}\label{pip}
	\Pi(p):=\sup _{\theta \in(0,1)} \lim _{m \rightarrow \infty} \sup _{q \geq 1} \mathbf{E}\left[\exp \left\{p \gamma \sup _{t \in[t_0,t_0+h]} \delta_\theta \bar{Y}_t^{(m, q)}\right\}\right]<\infty,
	\end{equation}
	where we use the following notations
\begin{equation}\label{thetaY}
\begin{split}
	 \delta_\theta Y^{(m, q)} & =\frac{\theta Y^{(m+q)}-Y^{(m)}}{1-\theta}, \quad \delta_\theta \tilde{Y}^{(m, q)}=\frac{\theta Y^{(m)}-Y^{(m+q)}}{1-\theta} \text { and} \\
		\delta_\theta \bar{Y}^{(m, q)} & =\left|\delta_\theta Y^{(m, q)}\right|+\left|\delta_\theta \tilde{Y}^{(m, q)}\right|.	
  \end{split} 
 \end{equation}
 The proof of uniform estimate for $\left(Y^{(m)}, Z^{(m)}\right)\in\mathbb{S} [t_0,t_0+h]\times \mathcal{H}^p [t_0,t_0+h]$ is postponed to the Appendix.

 Here we only prove that the total variation of $K^{(m)}$ is uniformly bounded. For each $m$,  consider the following of backward Skorokhod problem $\mathbb{BSP}^{r^{1,(m)}}_{l^{1,(m)}}(s^{1,(m)},a^{1,(m)})$ with two linear constraints. Here for any $(t,x)\in[t_0,t_0+h]\times \mathbb{R}$,
\begin{equation}
\begin{split}
&a^{1,(m)}:=\E[\eta],\\
&s^{1,(m)}_t:=\E
\left[\int_{t_0}^t f\left(s,Y_s^{(m-1)}, \mathbf{P}_{Y_s^{(m-1)}},Z^{(m)}_s\right)ds\right]=\E\left[y^{(m)}_{t_0}-y^{(m)}_t\right], \\ 
&l^{1,(m)}(t,x):=\E\left[L'\left(t,y^{(m)}_t-\E[y^{(m)}_t]+x\right)\right]=b_tx-p_t,\\
&r^{1,(m)}(t,x):=\E\left[R'\left(t,y^{(m)}_t-\E[y^{(m)}_t]+x\right)\right]=b_tx-q_t.
\end{split}
\end{equation}
By Remark \ref{BSPtoSP}, for any $t\in[t_0,t_0+h]$, we have $K^{1,(m)}_t=\bar{K}^{1,(m)}_h-\bar{K}^{1,(m)}_{t_0+h-t}$,  where $\bar{K}^{1,(m)}$ is the second component of the solution to the Skorokhod problem $\mathbb{SP}_{\bar{l}^{1,(m)}}^{\bar{r}^{1,(m)}}(\bar{s}^{1,(m)})$. Here, $\bar{s}^{1,(m)}$, $\bar{l}^{1,(m)}$ and $\bar{r}^{1,(m)}$ are given as follows 
\begin{equation}\label{barsbarlbarr}\begin{split}
&\bar{s}^{1,(m)}_{t}=a+s^{1,(m)}_{t_0+h}-s^{1,(m)}_{t_0+h-t}=\mathbf{E}[y^{(m)}_{t_0+h-t}],\\ 
&\bar{l}^{1,(m)}(t,x)=l^{1,(m)} (t_0+h-t,x), \ \bar{r}^{1,(m)}(t,x)=r^{1,(m)}(t_0+h-t,x).
\end{split}\end{equation}
We assume $\phi^{(m)},\ \psi^{(m)}$ are the unique solutions to the following equations, respectively:
$$\bar{l}^{1,(m)}(t,\bar{s}^{1,(m)}_t+x)=0, \ \bar{r}^{1,(m)}(t,\bar{s}^{1,(m)}_t+x)=0.$$
Therefore, we have 
$$\phi_{t}^{(m)}=\frac{p_{t}}{b_{t}}-\bar{s}_{t}^{1,(m)},\ \psi_{t}^{(m)}=\frac{q_{t}}{b_{t}}-\bar{s}_{t}^{1,(m)}. $$
According to Proposition $\ref{comparison}$,  we have 
$$\bar{K}_t^{1,(m),R} \geq \bar{K}_t^{(m),R},\  \bar{K}_t^{1,(m),L} \geq \bar{K}_t^{(m),L}, \ \forall t\in[t_0,t_0+h].$$
Combining with Proposition $\ref{totalvariation}$, we have
\begin{equation}\begin{split}
Var^{t_0+h}_{t_0}(K^{(m)})
&= Var^{h}_{0}(\bar{K}^{(m)})
=\bar{K}^{(m),R}_h+\bar{K}^{(m),L}_h\\
&\leq \bar{K}^{1,(m),R}_h+\bar{K}^{1,(m),L}_h\\
&\leq Var^{h}_{0}(\phi^{(m)})+Var^{h}_{0}(\psi^{(m)})\\
&\leq Var^{h}_{0}(\frac{p}{b})+Var^{h}_{0}(\frac{q}{b})+2Var^{h}_{0}(\bar{s}^{1,(m)}).\\
\end{split}   
\end{equation}
 For $Var^{h}_{0}(\bar{s}^{1,(m)})$, recalling that the definition of $\bar{s}^{1,(m)}$, we have 
\begin{equation}
    \begin{split}
    &Var^{h}_{0}(\bar{s}^{1,(m)})=Var_{t_0}^{t_0+h}(s^{1,(m)})\\
   & \leq \mathbf{E} \left[\int_{t_0}^{t_0+h}\left|f(s,Y_s^{(m-1)},\mathbf{P}_{Y_s^{(m-1)}},Z^{(m)}_s)\right|ds\right]\\
    &\leq \mathbf{E}\left[\int_{t_0}^{t_0+h}\alpha_s ds\right]+2\lambda \mathbf{E} \left[\int_{t_0}^{t_0+h}|Y^{(m-1)}_s|ds\right]+\frac{\gamma}{2}\mathbf{E} \left[\int_{t_0}^{t_0+h}|Z_s^{(m)}|^2 ds\right]< \infty.\\
    \end{split}
\end{equation}
Therefore, for each $m$, we have $Var^{t_0+h}_{t_0}(K^{(m)})< \infty$.\\
{\bf Step 3.} We show that $(Y^{(m)},Z^{(m)},K^{(m)})$ converges to $(Y,Z,K)\in \mathbb{S} [t_0,t_0+h]\times \mathcal{H}^p [t_0,t_0+h] \times BV[t_0,t_0+h]$. 
In fact, for any $p >1$ and $\theta \in(0,1)$, we have 
$$
\begin{aligned}
	& \limsup _{m \rightarrow \infty} \sup _{q \geq 1} \mathbf{E}\left[\sup _{t \in[t_0,t_0+h]}|Y_t^{(m+q)}-Y_t^{(m)}|^p\right] \\
    =&\limsup_{m \rightarrow \infty} \sup _{q \geq 1} \mathbf{E}\left[(1-\theta)^p \sup _{t \in[t_0,t_0+h]}\left|\frac{\theta Y_{t}^{(m+q)}-\theta Y_{t}^{(m+q)}+Y_{t}^{(m+q)}-Y_{t}^{(m)}}{1-\theta}\right|^p\right]\\
    \leq&  \limsup_{m \rightarrow \infty} \sup _{q \geq 1}2^{p-1}(1-\theta)^p \E\left[\left(\sup _{t \in[t_0,t_0+h]}\delta_\theta \bar{Y}^{(m, q)}_t\right)^p+\sup _{t \in[t_0,t_0+h]}|Y_{t}^{(m+q)}|^p\right]\\
    \leq &\limsup_{m \rightarrow \infty}\sup _{q \geq 1} 2^{p-1}(1-\theta)^p\left(\frac{p!}{\gamma^p}\E\left[\exp\left\{\gamma\sup _{t \in[t_0,t_0+h]}\delta_\theta \bar{Y}^{(m, q)}_t\right\}\right]+\sup _{m \geq 1} \mathbf{E}\left[\sup _{t \in[t_0,t_0+h]}|Y_t^{(m)}|^p\right]\right)\\
	\leq &2^{p-1}(1-\theta)^p\left(\frac{\Pi(1) p!}{\gamma^p}+\sup _{m \geq 1} \mathbf{E}\left[\sup _{t \in[t_0,t_0+h]}|Y_t^{(m)}|^p\right]\right),\\
\end{aligned}
$$
which together with Eqs. \eqref{uniformbounded}, \eqref{pip} and the arbitrariness of $\theta$ implies that
\begin{equation}\label{convergenceY}
\limsup _{m \rightarrow \infty} \sup _{q \geq 1} \mathbf{E}\left[\sup _{t \in[t_0,t_0+h]}|Y_t^{(m+q)}-Y_t^{(m)}|^p\right]=0,\ \forall p > 1.
\end{equation}

Applying Itô's formula to $|Y_t^{(m+q)}-Y_t^{(m)}|^2$, we have
$$
\begin{aligned}
	& \left[\int_{t_0}^{t_0+h}|Z_t^{(m+q)}-Z_t^{(m)}|^2 d t\right]^{\frac{p}{2}} \\
   \leq& C_p\left[\left(\int_{t_0}^{t_0+h}|Y_t^{(m+q)}-Y_t^{(m)}|^2 d t\right)^{\frac{p}{2}}+\left(\int_{t_0}^{t_0+h}|Y_t^{(m+q)}-Y_t^{(m)}||Z_t^{(m+q)}-Z_t^{(m)}|d B_t\right)^{\frac{p}{2}}+\left(\Delta_1^{(m, q)}\right)^{\frac{p}{2}}\right]\\
    \leq &C_p\left[\sup _{t \in[t_0,t_0+h]}|Y_t^{(m+q)}-Y_t^{(m)}|^{p} h^{{\frac{p}{2}}}+\sup _{t \in[t_0,t_0+h]}|Y_t^{(m+q)}-Y_t^{(m)}|^{{\frac{p}{2}}} |\Delta_2^{(m, q)}|^{{\frac{p}{2}}}\right.\\
    &\left. \qquad\qquad\qquad\qquad\qquad\qquad\qquad\qquad\qquad+\left(\int_{t_0}^{t_0+h}|Y_t^{(m+q)}-Y_t^{(m)}||Z_t^{(m+q)}-Z_t^{(m)}|d B_t\right)^{\frac{p}{2}}\right] \\
\end{aligned}
$$
where
\begin{align*}
  \Delta_1^{(m, q)}:&=\int_{t_0}^{t_0+h} 2\left|Y_t^{(m+q)}-Y_t^{(m)}\right| \left|f\left(t, Y_t^{(m+q-1)}, \P_{Y_t^{(m+q-1)}},Z_t^{(m+q)}\right) -f\left(t,Y_t^{(m-1)},\P_{Y_t^{(m-1)}}, Z_t^{(m)}\right)\right| dt\\
  &\qquad-\int_{t_0}^{t_0+h}\left|Y_t^{(m+q)}-Y_t^{(m)}\right|d\left(K_t^{(m+q)}-K_t^{(m)}\right),\\
  \Delta_2^{(m, q)}:&=\int_{t_0}^{t_0+h}2\left|f\left(t, Y_t^{(m+q-1)}, \P_{Y_t^{(m+q-1)}},Z_t^{(m+q)}\right)-f\left(t,Y_t^{(m-1)},\P_{Y_t^{(m-1)}}, Z_t^{(m)}\right)\right|d t\\
  & \quad +Var^{t_0}_{t_0+h}(K^{(m+q)})+Var^{t_0}_{t_0+h}(K^{(m)}).
\end{align*}

Taking expectations in both sides yields that
$$
\begin{aligned}
	& \mathbf{E}\left[\left(\int_{t_0}^{t_0+h}\left|Z_t^{(m+q)}-Z_t^{(m)}\right|^2 d t\right)^{\frac{p}{2}}\right] \\
   \leq & C_p h^{\frac{p}{2}} \E\left[\sup _{t \in[t_0,t_0+h]}\left|Y_t^{(m+q)}-Y_t^{(m)}\right|^p\right]+ C_p\mathbf{E}\left[\left|\Delta_2^{(m, q)}\right|^p\right]^{\frac{1}{2}} \mathbf{E}\left[\sup _{t \in[t_0,t_0+h]}\left|Y_t^{(m+q)}-Y_t^{(m)}\right|^p\right]^{\frac{1}{2}}\\
   &\quad+C_p\E\left[\int_{t_0}^{t_0+h}\left|Y_t^{(m+q)}-Y_t^{(m)}\right|^2\left|Z_t^{(m+q)}-Z_t^{(m)}\right|^2dt\right] \\
	\leq & C_p h^{\frac{p}{2}} \E\left[\sup _{t \in[t_0,t_0+h]}\left|Y_t^{(m+q)}-Y_t^{(m)}\right|^p\right]+ C_p\mathbf{E}\left[\left|\Delta_2^{(m, q)}\right|^p\right]^{\frac{1}{2}} \mathbf{E}\left[\sup _{t \in[t_0,t_0+h]}\left|Y_t^{(m+q)}-Y_t^{(m)}\right|^p\right]^{\frac{1}{2}}\\
   &\quad+C_p\E\left[\sup _{t \in[t_0,t_0+h]}\left|Y_t^{(m+q)}-Y_t^{(m)}\right|^{\frac{p}{2}}\left(\int_{t_0}^{t_0+h}\left|Z_t^{(m+q)}-Z_t^{(m)}\right|^2dt\right)^{\frac{p}{4}}\right] \\
   \leq & C_p h^{\frac{p}{2}} \E\left[\sup _{t \in[t_0,t_0+h]}\left|Y_t^{(m+q)}-Y_t^{(m)}\right|^p\right]+ C_p\mathbf{E}\left[\left|\Delta_2^{(m, q)}\right|^p\right]^{\frac{1}{2}} \mathbf{E}\left[\sup _{t \in[t_0,t_0+h]}\left|Y_t^{(m+q)}-Y_t^{(m)}\right|^p\right]^{\frac{1}{2}}\\
   &\quad+C_p\E\left[\frac{1}{2}\sup _{t \in[t_0,t_0+h]}\left|Y_t^{(m+q)}-Y_t^{(m)}\right|^{p}+\frac{1}{2}\left(\int_{t_0}^{t_0+h}\left|Z_t^{(m+q)}-Z_t^{(m)}\right|^2dt\right)^{\frac{p}{2}}\right] \\
\end{aligned}
$$
where we have used H\"{o}lder's inequality and $\text{B-D-G}$ inequality in the first inequality, Young inequality in the third inequality. $C_p$ is a constant depending on $p$ and varies from line to line.    Applying Eq. \eqref{convergenceY}, we obtain 
$$
\underset{m \rightarrow \infty}{\limsup } \sup _{q \geq 1} \mathbf{E}\left[\left(\int_{t_0}^{t_0+h}\left|Z_t^{(m+q)}-Z_t^{(m)}\right|^2 d t\right)^{\frac{p}{2}}\right]=0, \quad \forall p > 1
$$
Therefore, there exists a pair of processes $ (Y,Z) \in \mathbb{S} [t_0,t_0+h]\times \mathcal{H}^p [t_0,t_0+h]$ such that
\begin{equation}\label{YmZm}
 \lim _{m \rightarrow \infty} \mathbf{E}\left[\sup _{t \in[t_0, t_0+h]}\left|Y_t^{(m)}-Y_t\right|^p+\left(\int_{t_0}^{t_0+h}\left|Z_t^{(m)}-Z_t\right|^2 d t\right)^{\frac{p}{2}}\right]=0,\quad  \forall p > 1 .\end{equation}
Set
$$
K_t=Y_t-Y_{t_0}+\int_{t_0}^t f\left(s,Y_s,\P_{Y_s}, Z_s\right) d s-\int_{t_0}^t Z_s d B_s, \ \forall t\in[t_0,t_0+h].
$$
By Assumption $\ref{ass-4}$ and Eq. \eqref{YmZm} , we have
\begin{equation}\label{fmf}
\lim _{m \rightarrow \infty} \mathbf{E}\left[\int_{t_0}^{t_0+h}\left|f\left(t, Y_t^{(m-1)}, \mathbf{P}_{Y_t^{(m-1)}}, Z_t^{(m)}\right)-f\left(t, Y_t, \mathbf{P}_{Y_t}, Z_t\right)\right| d t\right]=0.  
\end{equation}
Consequently,
\begin{equation}\label{Km}
\lim_{m \rightarrow \infty}\mathbf{E}\left[\sup _{t \in[t_0,t_0+h]}\left|K_t-K_t^{(m)}\right|\right] = 0.
\end{equation}
In particular, we have $K_t:=\lim _{m \rightarrow \infty} K_t^{(m)}=\lim _{m \rightarrow \infty} \mathbf{E}\left[K_t^{(m)}\right]=\mathbf{E}\left[K_t\right]$ and then $K$ is a  continuous function with bounded variation.

{\bf Step 4.} We show that the limit processes $(Y,Z,K)$ are the solution to the quadratic $\text{MFBSDE}$ with double mean reflections \eqref{Yt0}. First, letting $m$ tend to infinity in the first two equations in \eqref{Ytm}, by Eqs. \eqref{YmZm} and \eqref{Km}, we obtain the first two equations in \eqref{Yt0}. It remains to show the Skorokhod condition holds. For this purpose, let $(Y',Z',K')$ be the solution to the doubly mean reflected BSDE on $[t_0,t_0+h]$ with terminal value $\eta$, constant coefficient $\{f(s,Y_s,\P_{Y_s},Z_s)\}_{s\in[t_0,t_0+h]}$ and loss functions $L,R$. It suffices to prove that $Y'\equiv Y$ and $K'\equiv K$. Let $(\widetilde{Y},\widetilde{Z})$ be the solution to the BSDE on $[t_0,t_0+h]$ with terminal value $\eta$ and constant coefficient $\{f(s,Y_s,\P_{Y_s},Z_s)\}_{s\in[t_0,t_0+h]}$. By the proof of  Proposition \ref{pro-1}, we have
\begin{equation}\label{Y't}
Y'_{t}=\widetilde{Y}_t+K'_{t_0+h}-K'_{t}=\eta+\int_{t}^{t_0+h}f(s,Y_s,\P_{Y_s},Z_s)ds-\int_{t}^{t_0+h}\widetilde{Z}_s dB_s+K'_{t_0+h}-K'_{t},    
\end{equation}
where $K'$ is the second component of the solution to $\mathbb{BSP}^{r}_{l}(s,a)$. 
Here, for any $(t,x)\in[t_0,t_0+h]\times \mathbb{R}$,
\begin{equation}
\begin{split}
&s_t:=\E\left[\int_{t_0}^t f(s,Y_s,\P_{Y_s},Z_s)ds\right], \ a:=\E[\eta],\\
&l(t,x):=\E\left[L(t,\widetilde{Y}_t-\E[\widetilde{Y}_t]+x)\right], \ r(t,x):=\E\left[R(t,\widetilde{Y}_t-\E[\widetilde{Y}_t]+x)\right].
\end{split}
\end{equation}
In view of Proposition \ref{continuity}, we have
$$
\begin{aligned}
 \sup_{t\in[t_0,t_0+h]}\left|K_t^{(m)}-K'_t\right|&=\sup_{t\in[t_0,t_0+h]}\left|\left(\bar{K}_h^{(m)}-\bar{K}_{t_0+h-t}^{(m)}\right)-\left(\bar{K}'_h-\bar{K}'_{t_0+h-t}\right)\right|\\
 &=\left|\bar{K}_h^{(m)}-\bar{K}'_h\right|+\sup_{t\in[t_0,t_0+h]}\left|\bar{K}_{t_0+h-t}^{(m)}-\bar{K}'_{t_0+h-t}\right|\\
 &=\left|\bar{K}_h^{(m)}-\bar{K}'_h\right|+\sup_{t\in[0,h]}\left|\bar{K}_{t}^{(m)}-\bar{K}'_{t}\right|\\
 &\leq 2\sup_{t\in[0,h]}\left|\bar{K}_{t}^{(m)}-\bar{K}'_{t}\right|\\
 &\leq \frac{2C}{c}\sup_{t\in[0,h]}\left|\bar{s}_t^{(m)}-\bar{s}_t\right|+\frac{2}{c}\left(\bar{L}_{h}\vee\bar{R}_{h}\right), \\     
\end{aligned}
$$
where $\bar{K}^{(m)}$ and $\bar{K}'$ are the second component of the solution to the Skorokhod problems $\mathbb{SP}_{\bar{l}^{(m)}}^{\bar{r}^{(m)}}(\bar{s}^{(m)})$ and $\mathbb{SP}_{\bar{l}}^{\bar{r}}(\bar{s})$.  Note that
\begin{align*}
\bar{s}_{t}&=\E\left[\eta\right]+\E\left[\int_{t_0+h-t}^{t_0+h}f(s,Y_s,\P_{Y_s},Z_s)ds\right],\\
\bar{s}_{t}^{(m)}&=\E\left[\eta\right]+\E\left[\int_{t_0+h-t}^{t_0+h}f(s,Y_s^{(m-1)},\P_{Y_s^{(m-1)}},Z^{(m)}_s ) ds\right].
\end{align*}
We have
$$\sup_{t\in[0,h]}|\bar{s}_{t}^{(m)}-\bar{s}_{t}|\leq \sup_{t\in[0,h]}\E\left[\int_{t_0+h-t}^{t_0+h}\left|f( s,Y_s^{(m-1)},\P_{Y_s^{(m-1)}},Z^{(m)}_s)-f(s,Y_s,\P_{Y_s},Z_s)\right|ds\right].$$
It follows from \eqref{fmf} that 
\begin{align*}
    \lim_{m\rightarrow \infty}\sup_{t\in[0,h]}|\bar{s}_{t}^{(m)}-\bar{s}_{t}|=0.
\end{align*}
Applying Assumption \ref{ass2}, we obtain that 
$$
\begin{aligned}
	\bar{L}_{h}&= \sup_{(t,x)\in[0,h]\times \mathbb{R}} \left|\mathbf{E}\left[L(t_0+h-t,y^{(m)}_{t_0+h-t}-\mathbf{E}[y^{(m)}_{t_0+h-t}]+x)\right]-\E\left[L(t_0+h-t,\widetilde{Y}_{t_0+h-t}-\mathbf{E}[\widetilde{Y}_{t_0+h-t}]+x)\right]\right|\\
	&\leq \sup_{t\in[0,h]}\mathbf{E}\left[C\left|y^{(m)}_{t_0+h-t}-\widetilde{Y}_{t_0+h-t}-\mathbf{E}\left[y^{(m)}_{t_0+h-t}-\widetilde{Y}_{t_0+h-t}\right]\right|\right]\\
	&\leq 2C\sup_{t\in[t_0,t_0+h]}\mathbf{E}\left[\left|y^{(m)}_{t}-\widetilde{Y}_t\right|\right].
\end{aligned}
$$
Recall that $(y^{(m)},z^{(m)})$ can be seen as the solution to the BSDE on $[t_0,t_0+h]$ with terminal value $\eta$ and constant coefficient $\{f(s,Y^{(m-1)}_s,\P_{Y^{(m-1)}_s},Z^{(m)}_s)\}_{s\in[t_0,t_0+h]}$. By the a priori estimates for BSDEs and using \eqref{YmZm}, we have 
\begin{align*}
    \lim_{m\rightarrow \infty}\bar{L}_{h}=0,
\end{align*}
 and similar result holds for $\bar{R}_{h}$. All the above analysis implies that 
\begin{align*}
\lim_{m\to\infty}\sup_{t\in[t_0,t_0+h]}\left|K_t^{(m)}-K'_t\right|= 0,
\end{align*}
which together with \eqref{Km} implies that $ K'\equiv K$. Replacing $K'$ by $K$ in \eqref{Y't} and comparing with the first equation in \eqref{Yt0}, by the uniqueness of BSDE, we obtain that $Y'\equiv Y$. Therefore, the Skorokhod condition holds. \\
{\bf Step 5.} For the general time interval $[0,T]$, similar as a standard $\text{BSDE}$ approach, we divide the time interval $[0,T]$ into finitely many small time intervals. On each small time interval, the previous analysis ensures the existence of a local solution. By stitching the local solutions, we get a global solution on the whole time interval. Uniqueness on the whole interval is a direct consequence of the uniqueness on each small interval.
The proof is complete. 
\end{proof}

 Next, we list an example of the sequence $\left\{K^{(m)}\right\}_{m=1}^{\infty}$ being uniformly bounded with respect to $m$. 
\begin{example}
For linear case, let $L(t,x)=x-4, \ R(t,x)=x+1$. Apparently, the sequence $\left\{K^{(m)}\right\}_{m=1}^{\infty}$ is uniformly bounded with respect to $m$. For nonlinear case,  let $L(t,x)=-\frac{1}{2}\frac{x^2}{1+\left|x\right|}+x-4, \ R(t,x)=\frac{1}{2}\frac{x^2}{1+\left|x\right|}+x+1 $, and  $L^{'}(t,x)=x-3, \ R^{'}(t,x)=x-1$. In this case, we can check that $L,R,L^{'},R^{'}$ satisfy Assumptions \ref{ass2} and \ref{assl'r'}, Thus, $\left\{K^{(m)}\right\}_{m=1}^{\infty}$ is uniformly bounded w.r.t. $m$.    
\end{example}

\section*{Acknowledgement}
    This work was supported  by the National Natural Science Foundation of China (No. 12301178), the Natural Science Foundation of Shandong Province for Excellent Young Scientists Fund Program (Overseas) (No. 2023HWYQ-049), the Fundamental Research Funds for the Central Universities and  the Qilu Young Scholars Program of Shandong University. 

\renewcommand\thesection{Appendix A}
\section{Construction by penalization in a special case}\label{append:A}
\renewcommand\thesection{A}

In this section, we apply a penalization method to construct the solution to the $\text{MFBSDE}$ with double mean reflections. However, due to the fact that the constraints only made for the distribution of $Y$ but not for the pointwise value of $Y$, we cannot make sure that the sequence of penalized MFBSDEs is monotone since the comparison theorem may not hold. For this reason, we only consider the following $\text{MFBSDE}$ with two linear reflections, whose parameters are denoted by $(\xi,f,l,r)$,
\begin{equation}\label{linearcase}
\begin{cases}
Y_t=\xi+\int_t^T f\left(s,Y_s,\P_{Y_s},Z_s,\P_{Z_s}\right)ds-\int_t^T Z_s dB_s+K_T-K_t, \\
l_t\leq \E\left[Y_t\right]\leq r_t, \\
K_t=K^l_t-K^r_t, \ K^l,K^r\in I[0,T]\textrm{ and } \int_0^T \left(\E[Y_t]-l_t\right)dK_t^l=\int_0^T \left(r_t-\E[Y_t]\right)dK^r_t=0.
\end{cases}
\end{equation}
Here, we assume that the obstacles $l,r$ satisfy the following assumption.

\begin{itemize}
\item[($H_{rl}$)] $r_t=\int_0^t a_s ds$, $l_t=\int_0^t b_s ds$ with $\int_0^T|a_t|^2 dt <\infty$, $\int_0^T|b_t|^2 dt<\infty$ and $l_t\leq r_t$, $t\in[0,T]$.
\end{itemize}


Now, consider the following penalized $\text{MFBSDEs}$:
\begin{equation}\begin{split}\label{panelization}
Y_t^n=&\xi+\int_t^T f(s,Y_s^n,\P_{Y_s^n},Z_s^n,\P_{Z^n_s})ds-\int_t^T Z_s^n dB_s\\
&+\int_t^T n(\E[Y_s^n]-l_s)^-ds-\int_t^T n(\E[Y_s^n]-r_s)^+ ds.
\end{split}\end{equation}
By Theorem 2.1 in \cite{LLZ}, the above equation admits a unique pair of solution $(Y^n,Z^n)\in \mathcal{S}^2\times \mathcal{H}^2$. We define $K^{n,l}_t:=\int_0^t n(\E[Y_s^n]-l_s)^- ds$, $K^{n,r}_t:=\int_0^t n(\E[Y_s^n]-r_s)^+ds$ and $K_t^n:=K^{n,l}_t-K^{n,r}_t$. We show that the solution to \eqref{linearcase} can be approximated by the solutions to \eqref{panelization}.

\begin{theorem}\label{existence and uniqueness}
Suppose that $f$ satisfies Assumption \ref{ass1} with $p=2$ and $r,l$ satisfy condition ($H_{lr}$). Given $\xi\in \mathcal{L}^2_T$ with $l_T\leq \E[\xi]\leq r_T$, the $\text{MFBSDE}$ with double mean reflections \eqref{linearcase} has a unique solution $(Y,Z,K)\in \mathcal{S}^2\times\mathcal{H}^2\times BV[0,T]$. Furthermore, $(Y,Z,K)$ is the limit of $(Y^n,Z^n,K^n)$.
\end{theorem}

Before proving Theorem \ref{existence and uniqueness}, we first present some a priori estimates similar with the classical reflected BSDEs.
\begin{proposition}\label{uniquenesslinear}
Suppose that $r,l$ satisfy ($H_{rl}$) and $f^i$ satisfy Assumption \ref{ass1} with $p=2$, $i=1,2$. Given $\xi^i\in \mathcal{L}^2_T$ with $l_T\leq \E[\xi]\leq r_T$, let $(Y^i,Z^i,K^i)$ be the solution to the doubly mean reflected $\text{MFBSDE}$ with parameters $(\xi^i,f^i,l,r)$. Then, there exists a constant $C$ depending on $\lambda,T$, such that 
\begin{align*}
\E\left[\int_0^T |\hat{Y}_t|^2dt\right]+\E\left[\int_0^T |\hat{Z}_t|^2 dt\right]\leq C\E\left[|\hat{\xi}|^2+\int_0^T |\hat{f}_t|^2dt\right],
\end{align*}
where $\hat{Y}_t=Y^1_t-Y^2_t$, $\hat{Z}_t=Z^1_t-Z^2_t$, $\hat{\xi}=\xi^1-\xi^2$ and $$\hat{f}_t=f^1(t,Y^1_t,\P_{Y^1_t},Z^1_t,\P_{Z^1_t})-f^2(t,Y^1_t,\P_{Y^1_t},Z_t^1,\P_{Z_t^1}).$$
\end{proposition}

\begin{proof}
 Let $\hat{K}_t=K^1_t-K^2_t$ and $\widehat{f^2_t}=f^2(t,Y^1_t,\P_{Y^1_t},Z^1_t,\P_{Z^1_t})-f^2(t,Y^2_t,\P_{Y^2_t},Z_t^2,\P_{Z^2_t})$. Applying It\^{o}'s formula to $e^{at}\hat{Y}_t^2$, where $a$ is a positive constant to be determined later, we have
\begin{equation}
\begin{split}\label{e1}
&|\hat{Y}_t|^2 e^{at}+\int_t^T ae^{as}|\hat{Y}_s|^2 ds+\int_t^T e^{as} |\hat{Z}_s|^2 ds\\
=&e^{aT}|\hat{\xi}|^2+\int_t^T 2e^{as}\hat{Y}_s \left(\hat{f}_s+\widehat{f^2_s}\right) ds+\int_t^T 2e^{as}\hat{Y}_s d\hat{K}_s-2\int_t^T e^{as}\hat{Y}_s\hat{Z}_sdB_s.
\end{split}
\end{equation}
By the assumption on $f^2$, the H\"{o}lder inequality and the property of 1-Wasserstein distance, we obtain that
\begin{equation}
\begin{split}\label{e2}
\int_t^T 2e^{as}\hat{Y}_s \widehat{f^2_s} ds\leq &\int_t^T 2\lambda e^{as}\left(|\hat{Y}_s|^2+|\hat{Y}_s\hat{Z}_s|+|\hat{Y}_s|W_1\left(\P_{Y^1_s},\P_{Y^2_s}\right)+|\hat{Y}_s|W_1\left(\P_{Z^1_s},\P_{Z^2_s}\right)\right) ds\\
\leq &\int_t^T 2\lambda e^{as}\left(|\hat{Y}_s|^2+|\hat{Y}_s\hat{Z}_s|+|\hat{Y}_s|\E[|\hat{Y}_s|]+|\hat{Y}_s|\E[|\hat{Z}_s|]\right) ds\\
\leq& \int_t^T (3\lambda+8\lambda^2)e^{as}|\hat{Y}_s|^2 ds+\int_t^T \lambda e^{as}\E[|\hat{Y}_s|^2]ds\\
&+\frac{1}{4}\int_t^T e^{as}|\hat{Z}_s|^2 ds+\frac{1}{4}\int_t^T e^{as}\E[|\hat{Z}_s|^2] ds,\\
\int_t^T 2e^{as}\hat{Y}_s \hat{f}_s ds\leq &\int_t^T e^{as}|\hat{Y}_s|^2 ds+\int_t^T e^{as}|\hat{f}_s|^2 ds.
\end{split}
\end{equation}
Set $a=\frac{3}{2}+4\lambda+8\lambda^2$. Plugging \eqref{e2} into \eqref{e1} and then taking expectations on both sides yield that 
\begin{align}\label{e3'}
\frac{1}{2}\E\left[\int_t^T e^{as}|\hat{Y}_s|^2 ds+\int_t^T e^{as} |\hat{Z}_s|^2 ds\right]\leq \E\left[e^{aT}|\hat{\xi}|^2+\int_t^T e^{as}|\hat{f}_s|^2 ds+ \int_t^T 2e^{as}\hat{Y}_s d\hat{K}_s\right].
\end{align}
It follows from the flat-off condition that 
\begin{equation}\label{e3}
\begin{split}
\E\left[\int_t^T e^{as}\hat{Y}_s d\hat{K}_s\right]&=\int_t^T e^{as}\big((\E[Y_s^1]-l_s)-(\E[Y_s^2]-l_s)\big)dK^{1,l}_s\\
&+\int_t^T e^{as}\big((r_s-\E[Y_s^1])-(r_s-\E[Y_s^2])\big)dK^{1,r}_s\\
&-\int_t^T e^{as}\big((\E[Y_s^1]-l_s)-(\E[Y_s^2]-l_s)\big)dK^{2,l}_s\\
&-\int_t^T e^{as}\big((r_s-\E[Y_s^1])-(r_s-\E[Y_s^2])\big)dK^{2,r}_s\leq 0.
\end{split}
\end{equation}
Combining Eqs. \eqref{e3'} and \eqref{e3}, we get the desired result.
\end{proof}

\begin{remark}
Note that the proof of Proposition \ref{uniquenesslinear} do not need $r,l$ are absolutely continuous. If $r,l\in C[0,T]$ satisfy $l_T\leq \E[\xi]\leq r_T$, the result of Proposition \ref{uniquenesslinear} still holds.
\end{remark}

In the following of this subsection, $C$ will always represent a positive constant  depending on $T,\lambda$, which may vary from line to line. 

\begin{proposition}\label{estimateYnZn}
There exists a constant $C$ independent of $n$, such that
\begin{align*}
&\sup_{t\in[0,T]}\E[\left|Y_t^n\right|^2]\leq C\left(\int_0^T a_s^2ds+\E\left[\xi^2\right]+\E\left[\int_0^T |f(s,0,\delta_0,0,\delta_0)|^2 ds\right]\right),\\
&\E\left[\int_0^T \left|Z_s^n\right|^2 ds\right]\leq C\left(\int_0^T a_s^2ds+\E\left[\xi^2\right]+\E\left[\int_0^T \left|f(s,0,\delta_0,0,\delta_0)\right|^2 ds\right]\right).
\end{align*}
\end{proposition}

\begin{proof}
Set $\widetilde{Y}^n_t=Y^n_t-r_t$. Applying It\^{o}'s formula to $e^{\beta t}|\widetilde{Y}^n_t|^2$, where $\beta$ is a positive constant to be determined later, we have
\begin{equation}\label{eq1.61}\begin{split}
&e^{\beta t}|\widetilde{Y}^n_t|^2+\int_t^T \beta e^{\beta s}|\widetilde{Y}^n_s|^2 ds+\int_t^T e^{\beta s}\left|Z_s^n\right|^2 ds\\
=&e^{\beta T}|\xi-r_T|^2+\int_t^T 2e^{\beta s}\widetilde{Y}^n_s\left(f(s,Y_s^n,\P_{Y^n_s},Z_s^n,\P_{Z^n_s})+a_s\right)ds-\int_t^T 2e^{\beta s}\widetilde{Y}_s^n Z_s^ndB_s\\
&+\int_t^T 2ne^{\beta s}\widetilde{Y}_s^n\left(\E[Y_s^n]-l_s\right)^-ds-\int_t^T 2ne^{\beta s}\widetilde{Y}_s^n(\E[\widetilde{Y}_s^n])^+ ds.
\end{split}\end{equation}
Simple calculation yields that 
\begin{equation}\label{eq1.62}\begin{split}
&2\widetilde{Y}^n_s\left(f(s,Y_s^n,\P_{Y^n_s},Z_s^n,\P_{Z^n_s})+a_s\right)\\
\leq &|f(s,r_s,\delta_{r_s},0,\delta_0)|^2+|a_s|^2+\frac{1}{4}|Z^n_s|^2+\frac{1}{4}W_1^2(\P_{Z^n_s},\delta_{0})+\lambda W_1^2(\P_{{Y}^n_s},\delta_{r_s})+(2+3\lambda+8\lambda^2)|\widetilde{Y}^n_s|^2\\
\leq &|a_s|^2+2|f(s,0,\delta_0,0,\delta_0)|^2+4\lambda^2|r_s|^2+\frac{1}{4}|Z^n_s|^2+\frac{1}{4}\E[|Z^n_s|^2]+\lambda \E[|\widetilde{Y}^n_s|^2]+(2+3\lambda+8\lambda^2)|\widetilde{Y}^n_s|^2.
\end{split}\end{equation}
Set $\beta=3+4\lambda+8\lambda^2$. Plugging Eq. \eqref{eq1.62} to Eq. \eqref{eq1.61} and taking expectations on both sides, noting that \begin{displaymath}
\E[\widetilde{Y}^n_s](\E[Y_s^n]-l_s)^-\leq 0,\ \E[\widetilde{Y}^n_s](\E[\widetilde{Y}_s^n])^+\geq 0,
\end{displaymath}
 we finally derive that
\begin{align*}
\E\left[|\widetilde{Y}_t^n|^2+\int_t^T |Z^n_s|^2 ds\right]\leq &C\E\left[|\xi-r_T|^2+\int_t^T |f(s,0,\delta_0,0,\delta_0)|^2ds+\int_t^T r_s^2ds+\int_t^T a_s^2ds\right]\\
\leq &C\E\left[|\xi|^2+\int_0^T |f(s,0,\delta_0,0,\delta_0)|^2ds+\int_0^T a_s^2ds\right].
\end{align*}
Recalling the definition of $\widetilde{Y}^n_t$, we obtain the desired result.
\end{proof}

In the following, we show that $K^{n,l}_T$ and $K^{n,r}_T$ are uniformly bounded.
\begin{proposition}\label{estimateYn-r}
There exists a constant $C$ independent of $n$, such that
\begin{align*}
n^2\int_0^T\left|(\E[Y_t^n]-r_t)^+\right|^2 dt\leq {C},\
n^2\int_0^T\left|(\E[Y_t^n]-l_t)^-\right|^2 dt\leq {C}.
\end{align*}
\end{proposition}

\begin{proof}
We only prove the first inequality since the second one can be proved in a same manner. We define $y^n_t:=\E[Y^n_t]$. Taking expectations on both sides of Eq. \eqref{panelization}, we have
\begin{align*}
y_t^n=\E[\xi]+\int_t^T \E\left[f\left(s,Y^n_s,\P_{Y^n_s},Z^n_s,\P_{Z^n_s}\right)\right]ds+\int_t^T n(y^n_s-l_s)^- ds-\int_t^T n(y^n_s-r_s)^+ds.
\end{align*} 
It is easy to check that 
\begin{align*}
d\left|(y^n_t-r_t)^+\right|^2=-2(y^n_t-r_t)^+ \left[ \E\left[f(t,Y^n_t,\P_{Y^n_t},Z^n_t,\P_{Z^n_t})\right]+a_t+n(y^n_t-l_t)^- -n(y^n_t-r_t)^+\right]dt.
\end{align*}
Noting that $y^n_T-r_T\leq 0$ and $(y^n_t-r_t)^+(y^n_t-l_t)^-\leq 0$, we have
\begin{align*}
&|(y^n_0-r_0)^+|^2+2n\int_0^T |(y^n_t-r_t)^+|^2dt\\
=&2\int_0^T (y^n_t-r_t)^+\left(\E\left[f(t,Y^n_t,\P_{Y^n_t},Z^n_t,\P_{Z^n_t})\right]+a_t\right) dt+2n\int_0^T (y^n_t-r_t)^+(y^n_t-l_t)^-dt\\
\leq &n\int_0^T  \left|(y^n_t-r_t)^+\right|^2dt+\frac{1}{n}\int_0^T \left(\E[f(t,Y^n_t,\P_{Y^n_t},Z^n_t,\P_{Z^n_t})]+a_t\right) ^2dt\\
\leq &n\int_0^T  \left|(y^n_t-r_t)^+\right|^2dt+\frac{C}{n}\int_0^T \left(\E[\left|f(t,0,\delta_0,0,\delta_0)\right|^2]+\E[\left|Y_t^n\right|^2]+\E[\left|Z_t^n\right|^2]+a^2_t\right)dt.
\end{align*}
Combining the above inequality and Proposition \ref{estimateYnZn} yields the desired result.
\end{proof}

Applying Proposition \ref{estimateYnZn}, \ref{estimateYn-r} and the B-D-G inequality, we obtain the following uniform estimate for $Y^n$.
\begin{proposition}\label{estimateKn}
There exists a constant $C$ independent of $n$, such that
\begin{align*}
\E\left[\sup_{t\in[0,T]}\left|Y_t^n\right|^2\right]\leq C.
\end{align*}
\end{proposition}


\begin{proof}[Proof of Theorem \ref{existence and uniqueness}]
Uniqueness is a direct consequence of Theorem \ref{th-1} with $p=2$. It remains to show the existence. We prove that there exists a triple of processes $(Y,Z,K)\in\mathcal{S}^2\times \mathcal{H}^2\times BV[0,T]$ such that
\begin{align}\label{equation6.25}
\E\left[\sup_{t\in[0,T]}\left|Y_t^n-Y_t\right|^2\right]+\E\left[\int_0^T \left|Z_s^n-Z_s\right|^2 ds\right]+\sup_{t\in[0,T]}\left|K^n_t-K_t\right|\rightarrow 0, \textrm{ as } n\rightarrow \infty.
\end{align} 
Besides, $(Y,Z,K)$ is the solution to \eqref{linearcase}.

For this purpose, for any positive integer $n$,  set $\hat{X}:=X^{n+1}-X^n$ for $X=Y,Z,K$. By a similar analysis as Eqs. \eqref{e1} and \eqref{e2}, choosing $a=\frac{1}{2}+4(\lambda+2\lambda^2)$, we have
\begin{align*}
&e^{at}|\hat{Y}_t|^2+\left(\lambda+\frac{1}{2}\right)\int_t^T e^{as} |\hat{Y}_s|^2ds +\frac{3}{4}\int_t^T e^{as}|\hat{Z}_s|^2ds\\
\leq &\lambda\int_t^T e^{as} \E[|\hat{Y}_s|^2]ds +\frac{1}{4}\int_t^T e^{as}\E[|\hat{Z}_s|^2]ds+2\int_t^T e^{as}\hat{Y}_sd\hat{K}_s-2\int_t^T e^{as}\hat{Y}_s\hat{Z}_sdB_s,
\end{align*}
Taking expectations on both sides yields that 
\begin{equation}\label{e30}
\sup_{t\in[0,T]}\E[|\hat{Y}_t|^2]+\E\left[\int_0^T (|\hat{Y}_s|^2+|\hat{Z}_s|^2)ds\right]\leq C\sup_{t\in[0,T]}\E\left[\int_t^T e^{as}\hat{Y}_sd\hat{K}_s\right].
\end{equation}

It is easy to check that for any $x,y\in\mathbb{R}$, we have
\begin{align*}
-nx^2+(2n+1)xy-(n+1)y^2=-n(x-\frac{2n+1}{2n}y)^2+\frac{y^2}{4n}.
\end{align*}
Set $y^n_t=\E[Y^n_t]$, $v^n_t=(y^n_t-l_t)^-$ and $u^n_t=(y_t^n-r_t)^+$. Then, we may calculate that
\begin{align*}
&[(y^{n+1}_s-l_s)-(y^n_s-l_s)]((n+1)v_s^{n+1}-nv_s^n)\\
\leq&-n|v_s^n|^2+(2n+1)v^n_s v^{n+1}_s-(n+1)|v_s^{n+1}|^2\leq \frac{|v^{n+1}_s|^2}{4n}
\end{align*}
and
\begin{align*}
&-[(y^{n+1}_s-r_s)-(y^n_s-r_s)]((n+1)u_s^{n+1}-nu_s^n)\\
\leq&-n|u_s^n|^2+(2n+1)u^n_s u^{n+1}_s-(n+1)|u_s^{n+1}|^2\leq \frac{|u^{n+1}_s|^2}{4n}.
\end{align*}
Then, we obtain that
\begin{equation}\label{eq1.81}\begin{split}
\E\left[\int_t^T e^{as}\hat{Y}_sd\hat{K}_s\right]=&\int_t^T e^{as}(y^{n+1}_s-y_s^n)\bigg(\big[(n+1)v_s^{n+1}-nv_s^n\big]-\big[(n+1)u_s^{n+1}-nu_s^n\big]\bigg)ds\\
 =&\int_t^T e^{as}\left[(y^{n+1}_s-l_s)-(y^n_s-l_s)\right]((n+1)v_s^{n+1}-nv_s^n)ds\\
 &-\int_t^T e^{as}\left[(y^{n+1}_s-r_s)-(y^n_s-r_s)\right]((n+1)u_s^{n+1}-nu_s^n)ds\\
 \leq &\frac{C}{n}\int_0^T \left(|v^{n+1}_s|^2+|u^{n+1}_s|^2\right)ds\leq \frac{C}{n^3},
\end{split}\end{equation}
where we have used Proposition \ref{estimateYn-r} in the last inequality. 
Plugging this estimate into \eqref{e30} indicates that
\begin{align}\label{eq1.84}
\sup_{t\in[0,T]}\E[|\hat{Y}_t|^2]+\E\left[\int_0^T \left(|\hat{Y}_s|^2+|\hat{Z}_s|^2\right)ds\right]\leq \frac{C}{n^3}.
\end{align}
Since $K^n$ is deterministic, it follows from \eqref{panelization} that
\begin{align*}
\hat{K}_T-\hat{K_t}=\E[\hat{Y}_t]-\E\left[\int_t^T \left(f\left(s,Y^{n+1}_s,\P_{Y^{n+1}_s},Z^{n+1}_s,\P_{Z^{n+1}_s}\right)-f\left(s,Y_s^n,\P_{Y^n_s},Z_s^n,\P_{Z^n_s}\right)\right)ds\right]
\end{align*}
and 
\begin{align}\label{eq1.84'}
\hat{Y}_t=\E_t\left[\int_t^T \left(f\left(s,Y^{n+1}_s,\P_{Y^{n+1}_s},Z^{n+1}_s,\P_{Z^{n+1}_s}\right)-f\left(s,Y_s^n,\P_{Y^n_s},Z_s^n,\P_{Z^n_s}\right)\right)ds\right]+\hat{K}_T-\hat{K_t}.
\end{align}
Using the H\"{o}lder inequality, the property of 1-Wasserstein distance and the estimate \eqref{eq1.84}, we obtain that
\begin{align*}
\sup_{t\in[0,T]}|\hat{K}_T-\hat{K_t}|^2\leq \frac{C}{n^3}.
\end{align*}
Noting the representation for $\hat{Y}$ in \eqref{eq1.84'}, by the B-D-G inequality, we have
\begin{align*}
\E\left[\sup_{t\in[0,T]}|\hat{Y}_t|^2\right]\leq \frac{C}{n^3}.
\end{align*} 
All the above analysis indicates that \eqref{equation6.25} holds.

It remains to prove the limit processes $(Y,Z,K)$ is indeed the solution to \eqref{linearcase}. The proof is similar with the proof of Theorem 5.5 in \cite{L}, so we omit it. The proof is complete.
\end{proof}

\begin{remark}
Compared with the penalty method for the reflected MFBSDE studied in \cite{CHM,DEH}, the generator $f$ in our case may depend on the distribution of both the solution $Y$ and $Z$.  Besides, we do not need any monotonicity condition for the generator $f$. 
\end{remark}

\renewcommand\thesection{Appendix B}
\section{Supplementary proofs}\label{append:B}
\renewcommand\thesection{B}

\begin{proof}[Proof of Eq. \eqref{uniformbounded}]
We first show that $$\sup_{m\geq 0}\mathbf{E}\left [ \exp\left \{ p\gamma\sup_{t \in[t_{0}, t_{0}+h]}\left | Y_t^{(m)} \right |  \right \}   \right ] <\infty.$$
According to  \eqref{re-2} , we have   
$$
\begin{aligned}
\sup_{t\in[t_0,t_0+h]}\left|Y^{(m)}_t\right|&=\sup_{t\in[t_0,t_0+h]}\left|y^{(m)}_t+\bar{K}^{(m)}_{t_0+h-t}\right|\\
  &=\sup_{t\in[t_0,t_0+h]}\left|y^{(m)}_t+(\bar{K}^{(m)}_{t_0+h-t}-\bar{K}^{3}_{t_0+h-t})+\bar{K}^{3}_{t_0+h-t}\right|\\
  &\leq\sup_{t\in[t_0,t_0+h]}\left|y^{(m)}_t\right|+\sup_{t\in[t_0,t_0+h]}\left|\bar{K}^{(m)}_{t_0+h-t}-\bar{K}^{3}_{t_0+h-t}\right|+\sup_{t\in[t_0,t_0+h]}\left|\bar{K}^{3}_{t_0+h-t}\right|\text {,}\\
\end{aligned}
$$
where $\bar{K}^3$ is the second component of the solution to the Skorokhod problem $\mathbb{SP}_{\bar{l}^3}^{\bar{r}^3}(\bar{s}^3)$. Here $\bar{s}^3=\mathbf{E}[\eta]$ and 
\begin{align*}
&\bar{l}^3(t,x)=l^3(t_0+h-t,x)=\mathbf{E}\left[L(t_0+h-t,x)\right], \\
&\bar{r}^3(t,x)=r^3(t_0+h-t,x)=\mathbf{E}\left[R(t_0+h-t,x)\right].
\end{align*}

In view of Proposition $\ref{continuity}$ , we have
$$
\begin{aligned}
\sup_{t\in[t_0,t_0+h]}\left|\bar{K}^{(m)}_{t_0+h-t}-\bar{K}^{3}_{t_0+h-t}\right|&=\sup_{t\in[0,h]}\left|\bar{K}^{(m)}_{t}-\bar{K}^{3}_{t}\right|\\
&\leq\frac{C}{c}\sup_{t\in[0,h]}\left|\bar{s}^{(m)}_{t}-\bar{s}^{3}_{t}\right|+\frac{1}{c}\left(\bar{L}^{3}_h\vee\bar{R}^{3}_h\right)\text{,}\\     
\end{aligned}
$$
where $\bar{L}^3_h=\sup_{(t,x)\in[0,h]\times\mathbb{R}}\left|\bar{l}^{(m)}(t,x)-\bar{l}^{3}(t,x)\right|$ and
$\bar{R}^3_h=\sup_{(t,x)\in[0,h]\times\mathbb{R}}\left|\bar{r}^{(m)}(t,x)-\bar{r}^{3}(t,x)\right|$.\\

Recalling \eqref{barsbarlbarr}, we have 
$$\sup_{t\in [0,h]}\left|\bar{s}^{(m)}_{t}-\bar{s}^{3}_{t}\right|\leq\sup_{t\in [0,h]}\left|\mathbf{E}[y^{(m)}_{t_0+h-t}]-\mathbf{E}[\eta]\right|\leq\sup_{t\in [t_0,t_0+h]}\mathbf{E}[|y^{(m)}_{t}|]+\mathbf{E}\left[\left|\eta\right|\right]$$
and
$$
\begin{aligned}
\bar{L}_h^3&=\sup_{(t,x)\in[0,h]\times\mathbb{R}}\left|\bar{l}^{(m)}(t,x)-\bar{l}^{3}(t,x)\right|\\
  &=\sup_{(t,x)\in[0,h]\times\mathbb{R}} \left|\mathbf{E}\left[L(t_0+h-t,y^{(m)}_{t_0+h-t}-\E[y^{(m)}_{t_0+h-t}]+x)\right]-\mathbf{E}\left[L(t_0+h-t,x)\right]\right|\\
 &\leq 2C\sup_{t\in [0,h]}\E[|y_{t_0+h-t}^{(m)}|]=2C\sup_{t\in [t_0,t_0+h]}\mathbf{E}[|y^{(m)}_{t}|].\\
\end{aligned}
$$
Similar estimate holds for $\bar{R}^{3}_h$. The above analysis imply that 
\begin{equation}\label{t_0+h-t}
\sup_{t\in[t_0,t_0+h]}\left|\bar{K}^{(m)}_{t_0+h-t}-\bar{K}^{3}_{t_0+h-t}\right|\leq 3\frac{C}{c}\sup_{t\in [t_0,t_0+h]}\mathbf{E}[|y^{(m)}_{t}|]+\frac{C}{c}\mathbf{E}\left[\left|\eta\right|\right].    
\end{equation}

Consequently, we have
\begin{equation}\label{supY}
    \sup_{t\in[t_0, t_0+h]}|Y^{(m)}_t|\leq\sup_{t\in[t_0, t_0+h]}|y^{(m)}_t|+3\frac{C}{c}\sup_{t\in[t_0, t_0+h]}\mathbf{E}[|y^{(m)}_t|]+B_2
\end{equation}
where $B_2:=\frac{C}{c}\mathbf{E}\left[\left|\eta\right|\right]+\sup_{t\in[0, T]}\left|\bar{K}^3_t\right|$. Similar to Eq. (24) and Eq. (25) in \cite{HMW},  we have 
$$\exp\left\{\gamma|y_{t}^{(m)}|\right\}\leq \mathbf{E}_{t}\left[\exp\left\{\gamma\left(\widetilde{\eta}+\lambda h\left(\sup_{t\in[t_0, t_0+h]}\left|Y^{(m-1)}_t\right|+\sup_{t\in[t_0, t_0+h]}\mathbf{E}\left[\left|Y^{(m-1)}_t\right|\right]\right)\right)\right\}\right]$$
and for any $m\geq 1, p > 1$, $t\in[t_0,t_0+h]$,
\begin{equation}\label{yineq}
\begin{split}
&\mathbf{E}\left[\exp\left\{p\gamma\sup_{t\in[t_0,t_0+h]}|y_{t}^{(m)}|\right\}\right]\\
\leq& 4\mathbf{E}\left[\exp\left\{p\gamma\left(\widetilde{\eta} +\lambda h\left(\sup_{t\in[t_0, t_0+h]}\left|Y^{(m-1)}_t\right|+\sup_{t\in[t_0, t_0+h]}\mathbf{E}\left[\left|Y^{(m-1)}_t\right|\right]\right)\right)\right\}\right] \\   
\leq& 4\mathbf{E}\left[\exp\left\{p\gamma\left(\widetilde{\eta}+\lambda h\sup_{t\in[t_0, t_0+h]}\left|Y^{(m-1)}_t\right|\right)\right\}\right]\mathbf{E}\left[\exp\left\{p\gamma\lambda h\sup_{t\in[t_0,t_0+h]}\left|Y^{(m-1)}_t\right|\right\}\right]\\
\leq& 4\mathbf{E}\left[\exp\left\{2p\gamma\widetilde{\eta}\right\}\right]^{\frac{1}{2}}\mathbf{E}\left[\exp\left\{2p \gamma \lambda h \sup_{t\in[t_0,t_0+h]}\left|Y^{(m-1)}_t\right|\right\}\right],\\
\end{split}  
\end{equation}
where $\widetilde{\eta}=|\eta |+\int_{t_0}^{t_0+h} \alpha_s ds $. Applying Lemma $\ref{le-2}$,  we obtain that
$$
\begin{aligned}
	 \mathbf{E}\left[\exp \left\{p \gamma \sup _{t\in[t_0,t_0+h]}\left|Y_t^{(m)}\right|\right\}\right] &\leq \widetilde{C}_2 \mathbf{E}\left[\exp \left\{\left(4+\frac{24C}{c}\right) p \gamma \lambda h \sup _{t\in[t_0,t_0+h]}\left|Y_t^{(m-1)}\right|\right\}\right] \\
 &\leq \widetilde{C}_2 \mathbf{E}\left[\exp \left\{ p \gamma  \sup _{t\in[t_0,t_0+h]}\left|Y_t^{(m-1)}\right|\right\}\right]^{\left(4+\frac{24C}{c}\right)\lambda h},
\end{aligned}
$$
where $$\widetilde{C}_2=4 \exp\left\{p \gamma B_2\right\} \mathbf{E}\left[\exp\left\{\left(4+\frac{24C}{c}\right)p\gamma\widetilde{\eta}\right\}\right]^{\frac{1}{2}}.$$
We set $\rho=\frac{1}{1-\left(4+\frac{24C}{c}\right)\beta h}$.  Then iterating the above procedure $m$ times yields that
$$\mathbf{E}\left[\exp \left\{p \gamma \sup _{t\in[t_0,t_0+h]}\left|Y_t^{(m)}\right|\right\}\right]\leq
\widetilde{C}_2^{\rho} \mathbf{E}\left[\exp\left\{24\frac{C}{c}p\gamma\left(\left|\eta\right|+\int_{t_0}^{t_0+h}\alpha_s ds\right)\right\}\right]^{\frac{\rho}{2}}<\infty.$$
It then follows that $$\sup_{m\geq 0}\mathbf{E}\left [ \exp\left \{ p\gamma\sup_{t\in[t_0,t_0+h]}\left | Y_t^{(m)} \right |  \right \}   \right ] <\infty.$$

Second, according to  inequality $(\ref{yineq})$, we have
$$\sup _{m \geq 0} \mathbf{E}\left[\exp \left\{p \gamma \sup _{t\in[t_0,t_0+h]}|y_t^{(m)}|\right\}\right]<\infty, \quad \forall p > 1.$$
Similar to \eqref{t_0+h-t} ,  we have 
$$
\begin{aligned}
\sup _{m \geq 0} \left|K_{t_0+h}^{(m)}\right| 
&=\sup _{m \geq 0} \left|\bar{K}_{h}^{(m)}-\bar{K}_{0}^{(m)}\right|\\
&=\sup _{m \geq 0} \left|\bar{K}_{h}^{(m)}-\bar{K}_{h}^{3}+\bar{K}_{h}^{3}\right|\\
&\leq 3\frac{C}{c}\sup_{m \geq 0}\mathbf{E}[\sup _{t\in[t_0,t_0+h]}|y_t^{(m)}|]+\frac{C}{c}\mathbf{E}\left[\left|\eta\right|\right]+\sup_{t\in[0,T]}\left|\bar{K}^3_t\right|<\infty \\
\end{aligned}
$$
Finally, since $Z^{(m)}=z^{(m)}$ (see Eq. \eqref{Zmzm}), applying  Corollary 4 in $\cite{BH08}$ to the quadratic $\text{BSDE}$ \eqref{ytm}, we obtain 
$$
\sup _{m \geq 0} \mathbf{E}\left[\left(\int_{t_0}^{t_0+h}\left|Z_t^{(m)}\right|^2 d t\right)^p\right]<\infty, \forall p > 1.
$$  
The proof is complete.
\end{proof}

\begin{proof}[Proof of Eq. \eqref{pip}]
 Actually, the proof is similar to the one for Eq. \eqref{bddexpdetltay}. In the following, we mainly highlight the differences. We first show that for any $m,q \geq 1$ , $p>1$ and $\theta\in\left(0,1\right)$
$$
\begin{aligned}
\mathbf{E} \left[\exp\left\{p\gamma\sup_{s\in[t_0,t_0+h]}\delta_{\theta}\bar{y}_s^{(m,q)}\right\}\right]&\leq4\mathbf{E}\left[\exp\left\{16p\gamma\zeta^{(m, q)} \right\}\right]^{\frac{1}{2}}\mathbf{E}\left[\exp\left\{16p\gamma\left(|\eta|+\chi^{(m, q)}\right)\right\}\right]^{\frac{1}{4}}\\
&\quad \times \mathbf{E}\left[\exp\left\{16p\gamma\lambda h\sup _{s \in[t_0, t_0+h]}\left|\delta_\theta Y_s^{(m-1, q)}\right|\right\}\right],   
\end{aligned}
$$
where the pair of processes $\left(\delta_\theta y^{(m, q)}, \delta_\theta z^{(m, q)}\right)$ are defined similarly as Eq. ($\ref{thetaY}$) and
$$
\begin{aligned}
	\zeta^{(m, q)} & =|\eta|+2\lambda h  \sup_m \mathbf{E}\left[\sup _{s \in[t_0, t_0+h]}\left|Y_s^{(m)}\right|\right]+\int_t^{t_0+h} \alpha_s d s+\lambda h\left(\sup _{s \in[t_0, t_0+h]}\left|Y_s^{(m-1)}\right|\right.\\
  &\left.\qquad+\sup _{s \in[t_0, t_0+h]}\left|Y_s^{(m+q-1)}\right|\right), \\
	\chi^{(m, q)} & =4 \lambda h\sup_m \mathbf{E}\left[\sup _{s \in[t_0, t_0+h]}\left|Y_s^{(m)}\right|\right] +\int_t^{t_0+h} \alpha_s d s+2 \lambda h\left(\sup _{s \in[t_0, t_0+h]}\left|Y_s^{(m-1)}\right|\right.\\
 &\left.\qquad+\sup _{s \in[t_0, t_0+h]}\left|Y_s^{(m+q-1)}\right|\right).
 \end{aligned}
 $$

 Recalling \eqref{ytm}, it is easy to check that $\left(\delta_\theta y^{(m, q)}, \delta_\theta z^{(m, q)}\right)$
satisfies the following $\text{BSDE}$
$$
\delta_\theta y_t^{(m, q)}=-\eta+\int_t^{t_0+h}\left(\delta_\theta f^{(m, q)}\left(s, \delta_\theta z_s^{(m, q)}\right)+\delta_\theta f_0^{(m, q)}(s)\right) d s-\int_t^{t_0+h} \delta_\theta z_s^{(m, q)} d B_s,
$$
where
$$
\begin{aligned} 
 \delta_\theta f_0^{(m, q)}(t)&= \frac{1}{1-\theta}\left(f\left(t, Y_t^{(m+q-1)}, \mathbf{P}_{Y_t^{(m+q-1)}}, z_t^{(m)}\right)-f\left(t, Y_t^{(m-1)}, \mathbf{P}_{Y_t^{(m-1)},}, z_t^{(m)}\right)\right), \\
\delta_\theta f^{(m, q)}(t, z)&=  \frac{1}{1-\theta}\left(\theta f\left(t, Y_t^{(m+q-1)}, \mathbf{P}_{Y_t^{(m+q-1)}}, z_t^{(m+q)}\right)\right. \\
& \left.\qquad \qquad \qquad \qquad \qquad-f\left(t, Y_t^{(m+q-1)}, \mathbf{P}_{Y_t^{(m+q-1)}},-(1-\theta) z+\theta z_t^{(m+q)}\right)\right) .\\
\end{aligned}
$$   
Similar to the proof of Lemma 9 in \cite{HMW}, we  obtain that for all $m,q \geq 1$ and $p>1$ 
$$
\begin{aligned}
& \mathbf{E} \left[\exp\left\{p\gamma\sup_{s\in[t_0,t_0+h]}\delta_{\theta}\bar{y}_s^{(m,q)}\right\}\right] \\ 
\leq& 4\mathbf{E}\left[\exp\left\{16p\gamma\zeta^{(m, q)} \right\}\right]^{\frac{1}{2}}\mathbf{E}\left[\exp\left\{8p\gamma\left(|\eta|+\chi^{(m, q)}+\lambda h\sup _{s \in[t_0, t_0+h]} \delta_\theta Y_s^{(m-1, q)}\right)\right\}\right]^{\frac{1}{2}}\\
& \times \mathbf{E}\left[\exp\left\{8p\gamma\lambda h\sup _{s \in[t_0, t_0+h]}\delta_\theta Y_s^{(m-1, q)}\right\}\right]^{\frac{1}{2}} \\
\leq &4\mathbf{E}\left[\exp\left\{16p\gamma\zeta^{(m, q)} \right\}\right]^{\frac{1}{2}}\mathbf{E}\left[\exp\left\{16p\gamma\left(|\eta|+\chi^{(m, q)}\right)\right\}\right]^{\frac{1}{4}}\mathbf{E}\left[\exp\left\{16p\gamma\lambda h\sup _{s \in[t_0, t_0+h]}\delta_\theta Y_s^{(m-1, q)}\right\}\right]^{\frac{1}{4}}\\
& \times \mathbf{E}\left[\exp\left\{16p\gamma\lambda h\sup _{s \in[t_0, t_0+h]}\delta_\theta Y_s^{(m-1, q)}\right\}\right]^{\frac{1}{2}} \\
\leq&4\mathbf{E}\left[\exp\left\{16p\gamma\zeta^{(m, q)} \right\}\right]^{\frac{1}{2}}\mathbf{E}\left[\exp\left\{16p\gamma\left(|\eta|+\chi^{(m, q)}\right)\right\}\right]^{\frac{1}{4}}\mathbf{E}\left[\exp\left\{16p\gamma\lambda h\sup _{s \in[t_0, t_0+h]}\delta_\theta Y_s^{(m-1, q)}\right\}\right].
\end{aligned}
$$

Secondly, through some elementary calculations,  we obtain that $$\delta_{\theta}\bar{Y}_{s}^{(m,q)}\leq \delta_{\theta}\bar{y}_{s}^{(m,q)}+6\frac{C}{c}\sup_{s\in[t_0,t_0+h]}\mathbf{E}[\delta_{\theta}\bar{y}_{s}^{(m,q)}]+2B_3, $$
where $B_3:=2\sup_{s\in[t_0,t_0+h]}\left|\bar{K}^3_s\right|+\left(\frac{C}{c}+1\right)\mathbf{E}\left[\left|\eta\right|\right]+12\frac{C}{c}\sup_{m}\mathbf{E}[\sup_{s\in[t_0,t_0+h]}|y_s^{(m)}|]<\infty$. Applying Lemma $\ref{le-2}$, we have 
$$
\begin{aligned}
 \mathbf{E}\left[\exp\left\{p\gamma\sup_{s\in[t_0,t_0+h]}\delta_{\theta}\bar{Y}_{s}^{(m,q)}\right\}\right]&\leq\widetilde{C}_3 \mathbf{E}\left[\exp\left\{2(1+6\frac{C}{c})16p\gamma\lambda h\sup_{s\in[t_0,t_0+h]}\delta_\theta \bar{Y}_s^{(m-1, q)}\right\}\right]\\
 & \leq\widetilde{C}_3\mathbf{E}\left[\exp\left\{p\gamma\sup_{s\in[t_0,t_0+h]}\delta_\theta \bar{Y}_s^{(m-1, q)}\right\}\right]^{(32+192\frac{C}{c})\lambda h},
\end{aligned}
$$
where $$\widetilde{C}_3=4\mathbf{E}\left[\exp\left\{16p\gamma\zeta^{(m, q)} \right\}\right]^{\frac{1}{2}} \exp\left\{2p\gamma B_3\right\}\mathbf{E}\left[16p\gamma(16+96\frac{C}{c})\left(\left|\eta\right|+\chi^{(m,q)}\right)\right]^{\frac{1}{4}}.$$ 
Set $\widetilde{\rho}=\frac{1}{1-(32+192\frac{C}{c})\lambda h}$. Iterating the above procedure $m$ times yields that
$$
\begin{aligned}
 &\mathbf{E}\left[\exp\left\{p\gamma\sup_{s\in[t_0,t_0+h]}\delta_{\theta}\bar{Y}_{s}^{(m,q)}\right\}\right]
 \leq\widetilde{C}_3^{\widetilde{\rho}} \mathbf{E}\left[\exp\left\{p\gamma\sup_{s\in[t_0,t_0+h]}\delta_\theta \bar{Y}_s^{(1, q)}\right\}\right]^{(32\lambda h+192\frac{C}{c}\lambda h)^{m-1}}.
\end{aligned}
$$
By the uniform estimate in  \eqref{uniformbounded}, we know that for any $\theta\in(0,1)$
$$\lim_{m\to\infty}\sup_{q\geq 1}\mathbf{E}\left[\exp\left\{p\gamma\sup_{s\in[t_0,t_0+h]}\delta_\theta \bar{Y}_s^{(1, q)}\right\}\right]^{(32\lambda h+192\frac{C}{c}\lambda h)^{m-1}}=1.$$
Recalling Eq. \eqref{uniformbounded} again, we obtain that 
$$
\begin{aligned}
\Pi(p)=\sup_{\theta\in(0,1)}\lim_{m\to\infty}\sup_{q\geq 1}\mathbf{E}\left[\exp\left\{p\gamma\sup_{s\in[t_0,t_0+h]}\delta_{\theta}\bar{Y}_{s}^{(m,q)}\right\}\right]\leq \sup_{m,q\geq 1} \widetilde{C}_3^{\widetilde{\rho}} < \infty.
\end{aligned}$$ 
\end{proof}

\end{document}